\documentclass[11pt]{article}%
\usepackage{amsmath}
\usepackage{amsfonts}
\usepackage{amssymb}
\usepackage{graphicx}%
\newtheorem{theorem}{Theorem}

\newtheorem{corollary}[theorem]{Corollary}

\newtheorem{definition}[theorem]{Definition}

\newtheorem{proposition}[theorem]{Proposition}
\newtheorem{remark}[theorem]{Remark}

\begin{document}

\begin{center}
{\LARGE Folding Difference and Differential Systems}

{\LARGE into Higher Order Equations}

\medskip

\centerline{H. SEDAGHAT \footnote{Department of Mathematics, Virginia Commonwealth University Richmond, Virginia, 23284-2014, USA; Email: h.sedaghat@discretedynamics.net}}
\end{center}

\vspace{2ex}

\begin{abstract}
A typical system of $k$ difference (or differential) equations can be
compressed, or folded into a difference (or ordinary differential) equation of
order $k$. Such foldings appear in control theory as the canonical forms of
the controllability matrices. They are also used in the classification of
systems of three nonlinear differential equations with chaotic flows by
examining the resulting jerk functions. The solutions of the higher order
equation yield one of the components of the system's $k$-dimensional orbits
and the remaining components are determined from a set of associated passive
equations. The folding algorithm uses a sequence of substitutions and
inversions along with index shifts (for difference equations) or higher
derivatives (for differential equations). For systems of two difference or
differential equations this compression process is short and in some cases
yields second-order equations that are simpler than the original system. For
all systems, the folding algorithm yields detailed amount of information about
the structure of the system and the interdependence of its variables. As with
two equations, some special cases where the derived higher order equation is
simpler to analyze than the original system are considered.

\end{abstract}

\bigskip

\section{Introduction}

It is common knowledge that a difference or ordinary differential equation of
order $k$ may be \textquotedblleft unfolded" in a standard way to a system of
$k$ first order difference or differential equations. A reverse process that
compresses, or \textquotedblleft folds" systems into higher order equations is
also possible, not just in rare instances. Folding linear systems in both
continuous and discrete time is seen in control theory; and folding nonlinear
systems of differential equations appears in the study of conditions that lead
to the occurrence of chaotic flows.

In control theory the \textquotedblleft controllability canonical form" is
precisely the folding, in the sense to be made precise here, of a
controllability matrix into a linear higher order equation, whether in
continuous or discrete time; see, e.g., \cite{Barn}, \cite{Elay}, \cite{Lasl}.
Using standard algebraic methods, a completely controllable system is found to
be equivalent to a linear equation (difference in discrete time, differential
in continuous time) whose order equals the rank of the controllability matrix.

In an entirely different line of research, in \cite{Eich} and \cite{Linz} a
variety of nonlinear differential systems displaying chaotic behavior are
studied and classified by converting them to ordinary differential equations
of order 3 that are defined by jerk functions (time rates of change of
acceleration). These systems include the well-known systems of Lorenz
\cite{lrnz} and R\H{o}ssler \cite{rslr} and most of the 19 autonomous,
minimally nonlinear systems introduced by Sprott in \cite{Sprt1} each
containing only a single quadratic term. These results define new categories
for distinguishing among a broad range of differential systems; for instance,
the systems of Lorenz and R\H{o}ssler can both be converted to third-order
ordinary differential equations using the same approach, but while
R\H{o}ssler's system folds globally with a jerk function having no
singularities, known jerk functions for Lorenz's system are not defined globally.

In this paper, we find that the aforementioned ideas in control theory and in
chaotic differential systems are special instances of the same concept,
namely, folding systems to equations. We take three further steps beyond the
aforementioned literature. First, we apply folding to systems generally,
whether linear or nonlinear, autonomous or not, in discrete time or in
continuous time. The basic idea behind our approach here is simple: starting
with a given system, a higher order equation is derived through a sequence of
\textit{substitutions}, \textit{inversions} and \textit{index shifts} (for
difference systems) or \textit{higher derivatives} (for differential ones).
These three actions constitute the main components of the \textit{folding
algorithm}. Although not the most elegant, an algorithmic approach is
preferred in the general context because folding a typical system with many
equations is an intricate and often tedious process.

For systems of two equations folding is a short, one-step process that is
practically useful in some cases where the resulting second-order equation is
simpler or more tractable than the original system. For a system of 3 or more
difference or differential equations, the folding algorithm, namely, the
iteration of the aforementioned 3-component folding process in principle
converts the system to a higher order equation with the inversion component
being the most technically uncertain part. Once a system is folded, results
formulated for the higher order equation may be used to analyze the solutions
of the system.

Taking a further step beyond the existing literature, we illustrate how
folding may be used in a precise way to examine the interdependence of
variables in a system. The aforementioned distinction between chaotic
differential systems of Lorenz, R\H{o}ssler and those in the Sprott family is
a case in point. In addition, non-fully controllable systems can be folded by
the same algorithm as the fully controllable ones and the folding process
indicates some of the differences between the two types of systems. More
generally, for parameter-dependent systems certain parameter values tend to
uncouple the system, at least partially. In the folding process, such
parameter values appear as singularity conditions and usually lead to
equations with orders that are lower than the number of equations in the
original system. Whether lower order foldings occur as a result of parameter
changes or are due to other structural features of a system, folding a system
provides new insights into the interdependence of variables.

Finally, in principle folding applies to non-autonomous systems in the same
way that it does to autonomous ones. Consideration of non-autonomous systems
is an important generalization on its own merit and it also leads to practical
benefits in some cases. We discuss examples of non-autonomous systems that
fold to higher order equations that are either autonomous or periodic, thus
making a greater range of methods applicable to those systems.

\section{Folding difference systems with two equations}

We begin with systems of two equations for which the folding process is
relatively simple. A (recursive or explicit) system of two first-order
difference equations is typically defined as%
\begin{equation}
\left\{
\begin{array}
[c]{c}%
x_{n+1}=f(n,x_{n},y_{n})\\
y_{n+1}=g(n,x_{n},y_{n})
\end{array}
\right.  \quad n=0,1,2,\ldots\label{s1}%
\end{equation}
where $f,g:\mathbb{N}\times D\rightarrow S$ are given functions, $\mathbb{N}$
is the set of non-negative integers, $S$ a nonempty set and $D\subset S\times
S$. An initial point $(x_{0},y_{0})\in D$ generates a (forward) orbit or
solution $\{(x_{n},y_{n})\}$ of (\ref{s1}) in the state-space $S\times S$
through the iteration of the function
\[
(n,x_{n},y_{n})\rightarrow(f(n,x_{n},y_{n}),g(n,x_{n},y_{n})):\mathbb{N}\times
D\rightarrow S\times S
\]
for as long as the points $(x_{n},y_{n})$ remain in $D.$ If (\ref{s1}) is
\textit{autonomous}, i.e., the functions $f,g$ do not depend on the index $n$
then $(x_{n},y_{n})=F^{n}(x_{0},y_{0})$ for every $n$ where $F^{n}$ denotes
the composition of the map $F(u,v)=(f(u,v),g(u,v))$ of $S\times S$ with itself
$n$ times.

A second-order, scalar difference equation in $S$ is defined as%
\begin{equation}
s_{n+2}=\phi(n,s_{n},s_{n+1}),\quad n=0,1,2,\ldots\label{e1}%
\end{equation}
where $\phi:\mathbb{N}\times D^{\prime}\rightarrow S$ is a given function and
$D^{\prime}\subset S\times S$. A pair of initial values $s_{0},s_{1}\in S$
generates a (forward) solution $\{s_{n}\}$ of (\ref{e1}) in $S$ if
$(s_{0},s_{1})\in D^{\prime}.$ As in the case of systems, if $\phi
(n,u,v)=\phi(u,v)$ is independent of $n$ then (\ref{e1}) is autonomous.

An equation of type (\ref{e1}) may be \textquotedblleft unfolded" to a system
of type (\ref{s1}) in a standard way; e.g.,%
\begin{equation}
\left\{
\begin{array}
[c]{l}%
s_{n+1}=t_{n}\\
t_{n+1}=\phi(n,s_{n},t_{n})
\end{array}
\right.  \label{s2}%
\end{equation}

In the system (\ref{s2}) the temporal delay in (\ref{e1}) is converted to an
additional variable in the state space. All solutions of (\ref{e1}) are
reproduced from the solutions of (\ref{s2}) in the form $(s_{n},s_{n+1}%
)=(s_{n},t_{n})$ so in this sense, higher order equations may be considered to
be special types of systems. In general, (\ref{e1}) has many possible
unfoldings into systems of two equations and (\ref{s2}) is not unique.

\subsection{An example\label{ilex}}

To highlight some issues of interest consider the system
\begin{equation}
\left\{
\begin{array}
[c]{l}%
x_{n+1}=x_{n}y_{n}\\
y_{n+1}=(a+bx_{n})/y_{n}%
\end{array}
\right.  \label{bs1}%
\end{equation}
of difference equations with $a,b\in\mathbb{R}$ (see Section \ref{inv} below
for how a system like (\ref{bs1}) may be derived). The domain $D$ of this
system is the complement of the x-axis ($y_{n}\not =0$) in the plane
$\mathbb{R}^{2}$. The first equation yields%
\begin{equation}
y_{n}=\frac{x_{n+1}}{x_{n}} \label{eyn0}%
\end{equation}
on the complement of the\textit{ y-axis} ($x_{n}\not =0$). Further, shifting
the index of $x_{n}$ in the first equation in (\ref{bs1}) by 1 yields
\begin{equation}
x_{n+2}=x_{n+1}y_{n+1}=x_{n+1}\left(  \frac{a+bx_{n}}{y_{n}}\right)
=x_{n+1}\left(  \frac{ax_{n}+bx_{n}^{2}}{x_{n+1}}\right)  =x_{n}(a+bx_{n})
\label{be1}%
\end{equation}

This\textit{ }second-order equation is effectively first-order in the sense
that its even and odd terms
\begin{align*}
x_{2k+2}  &  =x_{2k}(a+bx_{2k}),\quad x_{0}\text{ given and }k\geq0\\
x_{2k+1}  &  =x_{2k-1}(a+bx_{2k-1}),\quad x_{1}=x_{0}y_{0}\text{ given and
}k\geq1
\end{align*}
separately satisfy the same first-order recurrence $r_{n+1}=r_{n}(a+br_{n})$.
This recurrence is conjugate to the logistic equation when $a,b\not =0$ since%
\begin{equation}
r_{n+1}=ar_{n}\left(  1+\frac{b}{a}r_{n}\right)  \quad\text{so if }%
s_{n}=-\frac{b}{a}r_{n}\text{ then }s_{n+1}=as_{n}(1-s_{n}). \label{L}%
\end{equation}

This equation exhibits well-known dynamics in a bounded interval, including
chaotic behavior for a range of values of $b$.

The above equations determine the value of the x-component $x_{n}$ in the pair
$(x_{n},y_{n}).$ If $x_{n}\not =0$ for all $n$ then $y_{n}$\textit{ may be
calculated from (\ref{eyn0}) without having to use the second equation in the
system}. Otherwise, it is necessary to use the second equation of (\ref{bs1})
to calculate $y_{n}$ recursively. This situation occurs when $x_{0}=0$ in
which case, $x_{1}=x_{0}y_{0}=0$ also and we obtain the trivial solution of
(\ref{be1}) i.e., $x_{n}=0$ for all $n.$ In this case, the second equation of
the system gives%
\[
y_{n+1}=\frac{a+bx_{n}}{y_{n}}=\frac{a}{y_{n}}%
\]

All solutions $\{y_{0},a/y_{0},y_{0},a/y_{0},\ldots\}$ of this simple
first-order autonomous equation have period 2. For the planar system, these
translate into 2-cycles on the y-axis; i.e., the y-axis is invariant and all
solutions that start on it are (stable but non-attracting) 2-cycles. It is
worth mentioning that this situation does not occur in the exceptional case
$a=0.$

Next, we note that \textit{certain solutions of (\ref{be1}) do not yield a
solution for the system,} a situation that occurs because the domain of the
system is properly contained in the domain of (\ref{be1}), namely,
$\mathbb{R}^{2}$. If $a,b\not =0$, $x_{0}=-a/b$ and $y_{0}\not =0$ then from
(\ref{be1}) we find that $x_{2k}=0$ for all $k$ while $x_{2k-1}$ is determined
by the logistic equation (\ref{L}). This solution of (\ref{be1}) does not
generate a solution for (\ref{bs1}) since $y_{1}=0$.

The preceding analysis determines the \textit{global behavior} of the
solutions of the system (\ref{bs1}). In particular, it is easy to verify that
the results are consistent with information obtained through the standard
determination of fixed points and their local stability.

\subsection{Semi-invertiblility}

The preceding discussion in particular raises the following issues:

\begin{itemize}
\item It is possible to solve at least one equation of the system for one of
the variables. In the above example, the first equation is solved for $y_{n}$
as (\ref{eyn0}). This not only eliminates $y_{n}$ but also
makes it possible to determine $y_{n}$ from a solution of (\ref{be1}) or without having to solve the second equation of (\ref{bs1}) as a
(nonautonomous) difference equation;

\item To take advantage of the preceding item, certain domain and range
restrictions must be specified in order to pin down the precise relationship
between the solutions of the system and those of the higher order equation.
\end{itemize}

The following concept addresses the above two issues.

\begin{definition}
\label{smsl}Let $S,T$ be nonempty sets and consider a function $f:T\times
D\rightarrow S$ where $D\subset S\times S$. Then $f$ is
\textbf{semi-invertible} or \textbf{partially invertible} if there are sets
$M\subset D$, $M^{\prime}\subset S\times S$ and a function $h:T\times
M^{\prime}\rightarrow S$ such that
\begin{equation}
w=f(t,u,v)\Rightarrow v=h(t,u,w)\quad\text{for all }t\in T,\text{ }(u,v)\in
M\text{\ and }(u,w)\in M^{\prime}. \label{sol}%
\end{equation}

The function $h$ may be called a semi-inversion, or partial inversion of $f.$
If $f$ is independent of $t$ then $t$ is dropped from the above notation (see
the examples below).
\end{definition}

Semi-inversion refers more accurately to the \textit{solvability} of the
equation $w-f(t,u,v)=0$ for $v$. This recalls the implicit function theorem. A
general version that is based on the contraction principle holds in Banach
spaces \cite{L1} yields the following existence result when applied to the
function $w-f(u,v)$ (the $t$-independent case). The reader who is not
interested in this level of generality may simply think of the Banach space as
the real line.

\begin{proposition}
\label{ift}Let $S$ be a Banach space and $\mathcal{O}$ an open set in $S\times
S.$ Let $f:\mathcal{O}\rightarrow S$ be a $C^{k}$-differentiable function for
a positive integer $k$ and suppose that $f(u_{0},v_{0})=w_{0}$ for some
$(u_{0},v_{0})\in\mathcal{O}.$ If the partial derivative $\partial_{2}%
f(u_{0},v_{0})$ is invertible (as a linear map) then there exists $r>0$ and a
$C^{k}$ function $h_{0}:B_{r}(w_{0},u_{0})\rightarrow S$ defined on the open
ball of radius $r$ centered at $(w_{0},u_{0})\in S\times S$ such that
$f(u,h_{0}(w,u))=w$ for all $(w,u)\in B_{r}(w_{0},u_{0}).$ Further,%
\[
\text{if }\left\Vert (w-w_{0},u-u_{0})\right\Vert ,\left\Vert v-v_{0}%
\right\Vert <r\text{ and }w=f(u,v)\text{ then }v=h_{0}(w,u).
\]

\end{proposition}

The function $h_{0}$ in Proposition \ref{ift} is thus a semi-inversion of $f$
on the open ball $B_{r}(w_{0},u_{0})$ in the autonomous case. Often a function
$h$ that is defined on a larger set than the ball and equals $h_{0}$ on the
ball satisfies the above conditions, although the preceding result does not
imply that much in general.

Going in a different direction, a substantial class of semi-invertible
functions is supplied (globally) by the following idea. This notion is defined
in an algebraically general context although it is less generic analytically.

\begin{definition}
Let $(G,\ast)$ be a nontrivial group, $T$ a nonempty set and let $f:T\times
G\times G\rightarrow G$. If there are functions $f_{1},f_{2}:T\times
G\rightarrow G$ such that
\[
f(t,u,v)=f_{1}(t,u)\ast f_{2}(t,v)
\]

for all $u,v\in G$ and every $t\in T$ then we say that $f$ is
\textit{separable} on $G$ and write $f=f_{1}\ast f_{2}$ for short.
\end{definition}

Every affine function $f(n,u,v)=a_{n}u+b_{n}v+c_{n}$ with $a_{n},b_{n},c_{n}$
in a ring $R$ with identity is separable on the additive group $(R,+)$ for all
$n\geq1$ with $T=\mathbb{N}.$ Similarly, $f(t,u,v)=a(t)u+b(t)v+c(t)$ is
separable on $(\mathbb{R},+)$\ for all $t\in\mathbb{R}.$

Now, suppose that $f_{2}(t,\cdot)$ is a bijection for every $t\in T$ and
$f_{2}^{-1}(t,\cdot)$ is its inverse for each $t$; i.e., $f_{2}(t,f_{2}%
^{-1}(t,v))=v$ and $f_{2}^{-1}(t,f_{2}(t,v))=v$ for all $v.$ Evidently, a
separable function $f$ is semi-invertible if the component function
$f_{2}(t,\cdot)$ is a bijection for each fixed $t$, since for every $u,v,w\in
G$ and $t\in T$
\[
w=f_{1}(t,u)\ast f_{2}(t,v)\Rightarrow v=f_{2}^{-1}(t,[f_{1}(t,u)]^{-1}\ast
w)
\]
where map inversion and group inversion, both denoted by $-1,$ are
distinguishable from the context. In this case, an explicit expression for the
semi-inversion $h$ exists globally as%
\begin{equation}
h(t,u,w)=f_{2}^{-1}(t,[f_{1}(t,u)]^{-1}\ast w) \label{spsh}%
\end{equation}
with $M=M^{\prime}=G\times G.$ We summarize this observation as follows.

\begin{proposition}
\label{sep}Let $(G,\ast)$ be a nontrivial group and $f=f_{1}\ast f_{2}$ be
separable$.$ If $f_{2}(t,\cdot)$ is a bijection for each $t$ then $f$ is
semi-invertible on $G\times G$ with a semi-inversion uniquely defined by
(\ref{spsh}).
\end{proposition}

For the system (\ref{bs1}) in the last section, the function $f(u,v)=a+uv$ is
separable if $a=0$ on $G$ where $G$ is the group of nonzero real numbers under
ordinary multiplication. Here $f_{1}(u)=f_{2}(u)=u$ are both bijections,
making $f$ semi-invertible. If $a\not =0$ then $f$ is not separable but it is
semi-invertible with $h(u,w)=(w-a)/u$ on the complement of the y-axis in the
plane (the set $M^{\prime}$). Note that $f$ is semi-invertible for all values
of $a$ where $u\not =0.$

The function $f(u,v)=(v\circ u)^{-1}=u^{-1}\circ v^{-1}$ on a permutation
group $G$ is separable with $f_{2}(v)=v^{-1}$ a bijection. Thus, $f$ is
semi-invertible with $h(u,w)=(u\circ w)^{-1}.$ The function $f(u,v)=a_{n}%
u+b_{n}v+c_{n}$ with $a_{n},b_{n},c_{n}$ in a ring $R$ with identity is
separable on the additive group $(R,+)$ for all $n\geq1$ with $f_{1}%
(n,v)=a_{n}u+c_{n}$ and $f_{2}(n,v)=b_{n}v.$ If $b_{n}$ is a unit in $R$ for
all large $n$ then $f_{2}(n,\cdot)$ is a bijection for all such $n$ and thus,
$f$ is semi-invertible on $R$ with $h(n,u,w)=b_{n}^{-1}(w-a_{n}u-c_{n})$ for
large $n$. If $a_{n},b_{n}$ are not units for large $n$ then $f$ is separable
but not semi-invertible for either $u$ or $v$.

Next, suppose that $\{(x_{n},y_{n})\}$ is a solution of (\ref{s1}) in $D$.
Assume that one of the component functions in (\ref{s1}), say, $f$ is
semi-invertible. Then for $n\geq0$ there is a set $M\subset D$, a set
$M^{\prime}\subset S\times S$ and a function $h:\mathbb{N}\times M^{\prime
}\rightarrow S$ \ such that for $(x_{n},y_{n})\in M$ and $(x_{n},x_{n+1})\in
M^{\prime}$
\begin{equation}
x_{n+1}=f(n,x_{n},y_{n})\Rightarrow y_{n}=h(n,x_{n},x_{n+1}) \label{yn}%
\end{equation}

Therefore,
\begin{equation}
x_{n+2}=f(n+1,x_{n+1},y_{n+1})=f(n+1,x_{n+1},g(n,x_{n},y_{n}))=f(n+1,x_{n+1}%
,g(n,x_{n},h(n,x_{n},x_{n+1}))) \label{md}%
\end{equation}

For each $n\geq0$ the function%
\begin{equation}
\phi(n,u,w)=f(n+1,w,g(n,u,h(n,u,w))) \label{se2}%
\end{equation}
is defined on $\mathbb{N}\times M^{\prime}.$ If $\{s_{n}\}$ is the solution of
(\ref{e1}) with initial values $s_{0}=x_{0}$ and $s_{1}=x_{1}=f(0,x_{0}%
,y_{0})$ where $\phi$ defined by (\ref{se2}) then
\[
s_{2}=f(1,s_{1},g(0,s_{0},h(0,s_{0},s_{1})))=f(1,x_{1},g(0,x_{0}%
,h(0,x_{0},x_{1})))=f(1,x_{1},g(0,x_{0},y_{0}))=x_{2}%
\]

By induction, $s_{n}=x_{n}$ and by (\ref{yn}) $h(n,s_{n},s_{n+1}%
)=h(n,x_{n},x_{n+1})=y_{n}.$ It follows that
\begin{equation}
(x_{n},y_{n})=(s_{n},h(n,s_{n},s_{n+1})) \label{xy}%
\end{equation}
i.e., the solution $\{(x_{n},y_{n})\}$ of (\ref{s1}) can be obtained from a
solution $\{s_{n}\}$ of (\ref{e1}). These observations establish the following result.

\begin{theorem}
Suppose that $f$ in (\ref{s1}) is semi-invertible with $M,M^{\prime}$ and $h$
as in (\ref{sol}). Then each orbit of (\ref{s1}) in $D$ on which $h$ is
defined may be derived from a solution of (\ref{e1}) via (\ref{xy}) with
$\phi$ given by (\ref{se2}).
\end{theorem}

The following gives a name to the pair of equations that generate the
solutions of (\ref{s1}) in the above theorem.

\begin{definition}
\label{fld}(Folding) The pair of equations%
\begin{align*}
s_{n+2}  &  =\phi(n,s_{n},s_{n+1}),\quad s_{0}=x_{0},\ s_{1}=f(0,x_{0}%
,y_{0})\\
y_{n}  &  =h(n,x_{n},x_{n+1})
\end{align*}

where $\phi$ is defined by (\ref{se2}) is a folding of the system (\ref{s1}).
\end{definition}

Note that the equation for $y_{n}$ is \textit{passive} in the sense that it
simply evaluates a given function and no dynamics or iterations are involved.
Also observe that (\ref{s1}) may be considered an unfolding of the
second-order equation above. It is generally not equivalent to the standard
unfolding (\ref{s2}) of that equation.

If one of the component functions in the system is separable then a global
result is obtained using (\ref{spsh}).

\begin{corollary}
\label{sprd}Let $(G,\ast)$ be a nontrivial group and $f=f_{1}\ast f_{2}$ be
separable on $G\times G$. If $f_{2}(n,\cdot)$ is a bijection for every $n$
then every solution $\{(x_{n},y_{n})\}$ of (\ref{s1}) in $G$ is derived from a
solution $\{s_{n}\}$ of%
\begin{equation}
s_{n+2}=f_{1}(n+1,s_{n+1},g(n,s_{n},f_{2}^{-1}(n,[f_{1}(n,s_{n})]^{-1}\ast
s_{n+1})) \label{sde1}%
\end{equation}
that yields the x-component $x_{n}$ with the initial values $s_{0}%
=x_{0},\ s_{1}=f_{1}(0,x_{0})\ast f_{2}(0,y_{0}).$ Further, the solution
$\{s_{n}\}$ of (\ref{sde1}) yields, passively via (\ref{spsh}), the
y-component%
\begin{equation}
y_{n}=f_{2}^{-1}(n,[f_{1}(n,s_{n})]^{-1}\ast s_{n+1}). \label{ysp}%
\end{equation}

\end{corollary}

For example, if $f,g$ are as defined in (\ref{bs1}) then (\ref{sde1}) yields%
\[
s_{n+2}=s_{n+1}\frac{b+cs_{n}}{s_{n+1}/s_{n}}=s_{n}(b+cs_{n})
\]
which matches the earlier results above. In this case, $G$ is the group of
nonzero real numbers under multiplication.

The next result is, strictly speaking, a special case of Corollary \ref{sprd}.

\begin{corollary}
Let $a_{n},b_{n},c_{n}$ be sequences in a ring $R$ with identity and let
$g:\mathbb{N}\times R\times R\rightarrow R.$ If $b_{n}$ is a unit for all $n$
then the semilinear system
\begin{equation}
\left\{
\begin{array}
[c]{l}%
x_{n+1}=a_{n}x_{n}+b_{n}y_{n}+c_{n}\\
y_{n+1}=g(n,x_{n},y_{n})
\end{array}
\right.  \label{sls}%
\end{equation}
folds into the second-order difference equation%
\begin{align}
s_{n+2}  &  =\phi(n,s_{n},s_{n+1}),\quad\text{where: }s_{0}=x_{0}%
,\ s_{1}=a_{0}x_{0}+b_{0}y_{0}+c_{0},\label{phis}\\
\phi(n,u,w)  &  =c_{n+1}+a_{n+1}w+b_{n+1}g(n,u,b_{n}^{-1}(w-a_{n}%
u-c_{n}))\nonumber
\end{align}
For each solution $\{s_{n}\}$ of\ (\ref{phis}) the y-components of orbits of
(\ref{sls}) are given by the passive equation%
\[
y_{n}=b_{n}^{-1}(s_{n+1}-a_{n}s_{n}-c_{n}).
\]

\end{corollary}

For instance, consider the nonautonomous semilinear system%
\begin{equation}
\left\{
\begin{array}
[c]{l}%
x_{n+1}=(-1)^{n}x_{n}+y_{n}\\
y_{n+1}=\psi(x_{n})+(-1)^{n}y_{n}%
\end{array}
\right.  \label{nas}%
\end{equation}
on a ring $R$ with identity where $\psi:\mathbb{R}\rightarrow\mathbb{R}$ is an
arbitrary function (see Section \ref{inv} below for how a system like
(\ref{nas}) is derived). This system folds via (\ref{phis}) into the
(autonomous) equation
\begin{equation}
s_{n+2}=\psi(s_{n})-s_{n} \label{do2}%
\end{equation}
as is readily verified. For each initial point $(x_{0},y_{0})$ in
$\mathbb{R}^{2}$ the numbers $s_{0}=x_{0}$ and $s_{1}=x_{0}+y_{0}$ generate a
solution $\{s_{n}\}$ of the difference equation (\ref{do2}), which is similar
to (\ref{be1}) with respect to the odd- and even-numbered terms of its
solutions. The y-components of an orbit of the system (\ref{nas}) are then
calculated using the passive equation%
\[
y_{n}=s_{n+1}-(-1)^{n}s_{n}.
\]

\section{Folding differential systems with two equations}

Consider a planar system of two first-order differential equations of type%
\begin{equation}
\left\{
\begin{array}
[c]{c}%
dx/dt=f(t,x(t),y(t))\\
dy/dt=g(t,x(t),y(t))
\end{array}
\right.  \quad t\in\mathbb{R} \label{ds1}%
\end{equation}
where $f,g:\mathbb{R}\times D\rightarrow\mathbb{R}$ are given functions and
$D\subset\mathbb{R}\times\mathbb{R}.$ Although not considered here, the
functions $f,g$ may actually be defined on more general spaces where
differentiation is defined.

The folding process for the system (\ref{ds1}) of differential equations is
analogous to that discussed in the previous section for difference systems.
Suppose that $f$ is semi-invertible, so by (\ref{sol}) there are sets
$M,M^{\prime}$ and $h:\mathbb{R}\times M^{\prime}\rightarrow\mathbb{R}$ such
that
\begin{equation}
x^{\prime}(t)=f(t,x(t),y(t))\Rightarrow y(t)=h(t,x(t),x^{\prime}%
(t))\quad\text{for all }t. \label{des}%
\end{equation}

Under suitable differentiability hypotheses, using the subscript notation for
partial derivatives we obtain%
\begin{align*}
x^{\prime\prime}  &  =\frac{d}{dt}f(t,x,y)=f_{t}(t,x,y)+x^{\prime}%
f_{x}(t,x,y)+y^{\prime}f_{y}(t,x,y)\\
&  =f_{t}(t,x,y)+x^{\prime}f_{x}(t,x,y)+g(t,x,y)f_{y}(t,x,y)
\end{align*}

Now using (\ref{des}),%
\begin{equation}
x^{\prime\prime}=f_{t}(t,x,h(t,x,x^{\prime}))+x^{\prime}f_{x}%
(t,x,h(t,x,x^{\prime}))+g(t,x,h(t,x,x^{\prime}))f_{y}(t,x,h(t,x,x^{\prime})).
\label{de2}%
\end{equation}

Let $\phi(t,x,x^{prime})$ denote the right hand side of (\ref{de2}). From an
initial point $(x(t_{0}),y(t_{0}))\in D$ for a given $t_{0}\in\mathbb{R}$, the
ordinary differential equation%
\begin{equation}
s^{\prime\prime}=\phi(t,s,s^{\prime}) \label{de1}%
\end{equation}
with initial values $s(t_{0})=x(t_{0}),\ s^{\prime}(t_{0})=f(t_{0}%
,x(t_{0}),y(t_{0}))$ generates the x-component of the flow $(x(t),y(t))$ of
(\ref{ds1}). The y-component is calculated from the passive equation
\begin{equation}
y(t)=h(t,x(t),x^{\prime}(t)) \label{des0}%
\end{equation}

As in the discrete case, the pair of equations (\ref{de1}) and (\ref{des0})
constitute a folding of the differential system (\ref{ds1}).

In cases where (\ref{de2}) is simpler than (\ref{ds1}) a practical advantage
is gained. To illustrate, consider the non-autonomous system%
\begin{equation}
\left\{
\begin{array}
[c]{l}%
x^{\prime}=tx^{2}-y\\
y^{\prime}=ax+x^{2}+2t^{2}x^{3}-2txy
\end{array}
\right.  ,\quad a\in\mathbb{R} \label{des1}%
\end{equation}
(see Section \ref{inv} below for how a system like (\ref{des1}) is derived).
The first equation above may be solved for $y$ to give
\begin{equation}
y=tx^{2}-x^{\prime}=h(x,x^{\prime}). \label{des2}%
\end{equation}

Now (\ref{de2}) (or direct differentiation with respect to $t$ and
substitution) yields%
\[
x^{\prime\prime}=x^{2}+2txx^{\prime}-(2t^{2}x^{3}+x^{2}+ax-2txy)=x^{2}%
+2txx^{\prime}-2t^{2}x^{3}-x^{2}-ax+2tx(tx^{2}-x^{\prime})=-ax
\]

The linear autonomous equation $x^{\prime\prime}+ax=0$ is elementary and its
solution (depending on the parameter value $a$) determines $x(t).$ Then $y(t)$
is calculated from the passive equation (\ref{des2}).

%\begin{remark}
The role of $a$ as a bifurcation parameter is not obvious from (\ref{des1})
directly. Further, the important issue of the \textit{existence of solutions}
for a system of differential equations is sometimes more easily addressed when
(\ref{de2}) is a simple equation, as in the above case.

We close this section by pointing out the following relevant facts.

\begin{enumerate}
\item The above discussion may be generalized to real-valued functions $f,g$
on Banach spaces with $x$ and $y$ both from $\mathbb{R}$ into a given Banach
space. Differentiability in Banach spaces is discussed in \cite{L1} so it is
not necessary to derive the technical details. However, we limit our work on
differential systems here to the real case.

\item The proper discrete analog for the differential system (\ref{ds1}) is
not (\ref{s1}) but rather, the difference system%
\begin{equation}
\left\{
\begin{array}
[c]{c}%
\Delta x_{n}=f(n,x_{n},y_{n})\\
\Delta y_{n}=g(n,x_{n},y_{n})
\end{array}
\right.  \label{del2}%
\end{equation}
where $\Delta x_{n}=x_{n+1}-x_{n}$ and similarly for $\Delta y_{n}.$ Although
(\ref{del2}) can be cast in the form (\ref{s1}) by slightly modifying $f,g$ it
is useful to keep the distinction between (\ref{s1}) and (\ref{del2}) in mind
when comparing difference and differential systems; see Remark \ref{tok} below.
\end{enumerate}

\section{Folding linear systems into equations}

The procedure described in the preceding sections applies, in particular, to
linear systems. Where there is no loss of clarity, we use the term
\textquotedblleft linear" loosely to include nonhomogeneous or affine
equations. Both of the coordinate functions in linear systems are separable
over the additive group of an underlying ring, though the requirement that one
of the constituent functions be a bijection is a multiplicative issue. Thus
semi-inversion requires the full ring structure.

\subsection{Linear difference systems and equations}

Consider the linear (nonhomogeneous, nonautonomous) system%
\begin{equation}
\left\{
\begin{array}
[c]{c}%
x_{n+1}=a_{n}x_{n}+b_{n}y_{n}+\alpha_{n}\\
y_{n+1}=c_{n}x_{n}+d_{n}y_{n}+\beta_{n}%
\end{array}
\right.  \label{dl}%
\end{equation}
where $a_{n},b_{n},c_{n},d_{n},\alpha_{n},\beta_{n}$ are given sequences in a
ring $R$ with identity$.$ Then shifting the index by 1 in the first equation
above gives%
\begin{align}
x_{n+2}  &  =a_{n+1}x_{n+1}+b_{n+1}y_{n+1}+\alpha_{n+1}\label{lin1}\\
&  =a_{n+1}x_{n+1}+b_{n+1}c_{n}x_{n}+b_{n+1}d_{n}y_{n}+b_{n+1}\beta_{n}%
+\alpha_{n+1}\nonumber
\end{align}

If $b_{n}$ is a \textit{unit }in $R$ for each $n$ then $f$ is semi-invertible
and%
\[
w=f(n,u,v)=a_{n}u+b_{n}v+\alpha_{n}\Rightarrow v=b_{n}^{-1}(w-a_{n}%
u-\alpha_{n}).
\]

Thus the y-component $y_{n}$ is evaluated as
\begin{equation}
y_{n}=b_{n}^{-1}(x_{n+1}-a_{n}x_{n}-\alpha_{n})=h(n,x_{n},x_{n+1}).
\label{linyn}%
\end{equation}

Substituting this into (\ref{lin1}) and rearranging terms yields%
\begin{equation}
x_{n+2}=A_{n}x_{n+1}+B_{n}x_{n}+C_{n} \label{aln}%
\end{equation}
where%
\begin{align*}
A_{n}  &  =a_{n+1}+b_{n+1}d_{n}b_{n}^{-1},\\
B_{n}  &  =b_{n+1}(c_{n}-d_{n}b_{n}^{-1}a_{n}),\\
C_{n}  &  =b_{n+1}(\beta_{n}-d_{n}b_{n}^{-1}\alpha_{n})+\alpha_{n+1}.
\end{align*}

Note that the sequences $\alpha_{n},\beta_{n}$ do not appear as bound
coefficients in (\ref{aln}) but within the free term $C_{n}$. If $\{x_{n}\}$
is a solution of this difference equation in $R$ then $y_{n}$ is passively
calculated from (\ref{linyn}).

It is clear that if $c_{n}$ is a unit for every $n$ (rather than $b_{n}$) then
we may switch the roles of x and y in the preceding discussion.

If neither $b_{n}$ nor $c_{n}$ are units for infinitely many $n$ then
(\ref{linyn}), or its analog for $c_{n}$, are not available so we proceed as
follows: the second equation of (\ref{dl}) yields%
\begin{equation}
b_{n+1}y_{n+1}=b_{n+1}c_{n}x_{n}+b_{n+1}d_{n}y_{n}+b_{n+1}\beta_{n}
\label{tmp0}%
\end{equation}
which may be used in (\ref{lin1}) to obtain%
\begin{equation}
x_{n+2}=a_{n+1}x_{n+1}+b_{n+1}c_{n}x_{n}+b_{n+1}d_{n}y_{n}+b_{n+1}\beta
_{n}+\alpha_{n+1} \label{tmp}%
\end{equation}

Next, multiply the first equation in (\ref{dl}) by $d_{n}$ to obtain%
\[
d_{n}x_{n+1}=d_{n}a_{n}x_{n}+d_{n}b_{n}y_{n}+d_{n}\alpha_{n}%
\]

Finally, \textit{if }$b_{n}$\textit{ is central} (commutes with all elements
of $R$) for every $n$ or if $R$ is a commutative ring then the above equality
yields%
\[
b_{n}d_{n}y_{n}=d_{n}x_{n+1}-d_{n}a_{n}x_{n}-d_{n}\alpha_{n}%
\]
and this together with (\ref{tmp}) produces%
\begin{align}
b_{n}x_{n+2}  &  =b_{n}a_{n+1}x_{n+1}+b_{n}b_{n+1}c_{n}x_{n}+b_{n+1}%
d_{n}x_{n+1}-b_{n+1}d_{n}a_{n}x_{n}-b_{n+1}d_{n}\alpha_{n}+b_{n}b_{n+1}%
\beta_{n}+b_{n}\alpha_{n+1}\nonumber\\
&  =(b_{n}a_{n+1}+b_{n+1}d_{n})x_{n+1}+b_{n+1}(b_{n}c_{n}-d_{n}a_{n}%
)x_{n}-b_{n+1}(d_{n}\alpha_{n}-b_{n}\beta_{n})+b_{n}\alpha_{n+1} \label{bnu1}%
\end{align}

If $\{x_{n}\}$ is a solution of (\ref{bnu1}) then the y-components are
obtained from the second equation of (\ref{dl}) as%
\begin{equation}
y_{n+1}=d_{n}y_{n}+(c_{n}x_{n}+\beta_{n}) \label{o1n0}%
\end{equation}

Though not as easy to use as (\ref{linyn}) because (\ref{o1n0}) is not a
passive equation, it is a linear and first-order difference equation that is
not difficult to solve.

If $b_{n}=b$ is a constant (albeit not a unit) then (\ref{bnu1}) takes a
simpler form in which the term $x_{n+2}$ is free. From the first equation of
(\ref{dl}) and (\ref{lin1}) we obtain%
\begin{align}
x_{n+2}  &  =a_{n+1}x_{n+1}+by_{n+1}+\alpha_{n+1}\nonumber\\
&  =a_{n+1}x_{n+1}+bc_{n}x_{n}+d_{n}by_{n}+b\beta_{n}+\alpha_{n+1}\nonumber\\
&  =a_{n+1}x_{n+1}+bc_{n}x_{n}+d_{n}(x_{n+1}-a_{n}x_{n}-\alpha_{n})+b\beta
_{n}+\alpha_{n+1}\nonumber\\
&  =(a_{n+1}+d_{n})x_{n+1}-(d_{n}a_{n}-bc_{n})x_{n}-d_{n}\alpha_{n}+b\beta
_{n}+\alpha_{n+1} \label{rcom}%
\end{align}

Once a solution $\{x_{n}\}$ of the above equation is obtained it is necessary
to solve the first-order equation (\ref{o1n0}) for $y_{n}.$ Unlike
(\ref{bnu1}), Eq. (\ref{rcom}) is recursive and its solutions may be obtained
by iteration.

The next result summarizes the preceding discussion.

\begin{proposition}
Consider the linear difference system (\ref{dl}) over a commutative ring with identity.

(a) If $b_{n}$ is a unit for all $n$ then the linear difference equation
(\ref{aln}) together with the passive equation (\ref{linyn}) constitute a
folding of (\ref{dl}).

(b) If $b_{n}$ is not a unit for infinitely many $n$ then (\ref{dl}) folds
into the (non-recursive) linear difference equation (\ref{bnu1}) that yields
the x-component of its orbits while its y-component is determined as a
solution of (\ref{o1n0}).

(c) If $b_{n}=b$ is a constant for all $n$ then (\ref{dl}) folds into the
recursive linear equation (\ref{rcom}) that yields the x-component of its
orbits while its y-component is determined either (i) passively from
(\ref{linyn}) if $b$ is a unit or (ii) as a solution of (\ref{o1n0}) if $b$ is
not a unit.
\end{proposition}

\subsection{Linear difference systems with periodic foldings}

We now consider linear systems that fold into equations with periodic
coefficients. \textit{We do not assume that the parameters of the linear
system are periodic}; it is only necessary that the derived second-order
equation has periodic coefficients so the results in this section apply to
many linear systems with non-periodic parameters as well. We consider only
difference systems in this section. Differential systems that fold into a
second-order ordinary differential equation with periodic coefficients may be
studied with the aid of the continuous-time Floquet theory \cite{Hln}.

Suppose that the coefficients $a_{n},b_{n},c_{n},d_{n}$\ in the system
(\ref{dl})\ are all sequences where $b_{n}$ is a unit for each $n$ (to limit
the range of possible cases) and the coefficients $A_{n},B_{n}$ in (\ref{aln})
are periodic with period $p$ i.e., $A_{n+p}=A_{n},\ B_{n+p}=B_{n}$. In
particular this is true if all the coefficients in the system are periodic
with period $p$ (not necessarily prime or minimal). However, non-periodic
systems may also fold into periodic equations of type (\ref{aln}). In fact, if
$A_{n},B_{n}$ are sequences with period $p$ and $b_{n},d_{n}$ are arbitrary
sequences then defining
\[
a_{n}=A_{n-1}-b_{n}d_{n-1}b_{n-1}^{-1},\quad c_{n}=b_{n+1}^{-1}B_{n}%
+d_{n}b_{n}^{-1}a_{n}%
\]
ensures that the homogeneous part of (\ref{aln}) has periodic coefficients.
Other combinations of system parameters that yield periodic $A_{n}$ and
$B_{n}$ are possible.

To study the solutions of (\ref{aln}) with periodic coefficients, one possible
approach is to unfold it back to a system and then use the Floquet theory
adapted to the discrete case; see, e.g., \cite{DHln}, \cite{SvD}.
Alternatively, we may use the eigensequence method in \cite{SL} which applies
direcly to (\ref{aln}) without the need for unfolding it and further, it works
whether $C_{n}$ is periodic or not. We need only find an eigensequence for the
homogeneous part of (\ref{aln}) in the ring $R$ (in particular, an
eigensequence with period $p$). Any eigensequence, periodic or not, yields a
semiconjugate factorization \cite{fsor} of the second-order equation into a
pair of first-order ones.

An eigensequence of period $p$ exists in $R$ if the characteristic equation of
(\ref{aln}), i.e., the first-order quadratic difference equation%
\begin{equation}
r_{n+1}r_{n}=A_{n-1}r_{n}+B_{n-1} \label{eper}%
\end{equation}
has a solution of period $p$ in the ring $R$ for some initial value $r_{1}\in
R$. The following result is from \cite{SL}.

\begin{theorem}
\label{per}Let $R$ be a ring with identity 1 and for $j=1,2,\ldots,p$, let
$\alpha_{j},\beta_{j}$ be obtained by iteration from (\ref{aln}) subject to
the initial values%
\[
\alpha_{0}=0,\ \alpha_{1}=1;\quad\beta_{0}=1,\ \beta_{1}=0.
\]
If a root $r_{1}$ of the quadratic polynomial%
\begin{equation}
r\alpha_{p}r+r\beta_{p}-\alpha_{p+1}r-\beta_{p+1}=0 \label{p2q1}%
\end{equation}
is a unit in $R$ and the recurrence%
\begin{equation}
r_{j+1}=A_{j-1}+B_{j-1}r_{j}^{-1} \label{epr}%
\end{equation}
\noindent also generates units $r_{2},\ldots,r_{p}$ in $R$ then $\{r_{n}%
\}_{n=1}^{\infty}$ is a unitary eigensequence of (\ref{aln}) with preiod $p$
that yields the triangular system of first-order equations (a semiconjugate
factorization)%
\begin{align*}
t_{n+1}  &  =C_{n-1}-B_{n-1}r_{n}^{-1}t_{n},\quad t_{1}=x_{1}-r_{1}x_{0}\\
x_{n+1}  &  =r_{n+1}x_{n}+t_{n+1}.
\end{align*}

\end{theorem}

The polynomial in (\ref{p2q1}) simplifies further if the coefficients
$A_{j},B_{j}$ are in the center of $R$. Then $\alpha_{j},\beta_{j}$ are also
in the center of $R$ so (\ref{p2q1}) reduces to
\begin{equation}
\alpha_{p}r^{2}+(\beta_{p}-\alpha_{p+1})r-\beta_{p+1}=0. \label{p2qc}%
\end{equation}
\ 

If $A_{n}=A$ and $B_{n}=B$ are constants then $p=1$ and the quadratic
polynomial (\ref{p2qc}) further reduces to $r^{2}-Ar-B=0$. This is
recognizable as the characteristic polynomial of the autonomous linear
equation of order 2.

For illustration suppose that the parameters of (\ref{dl}) are real with
$b_{n}\not =0$ for all $n$ and the following equalities hold%
\begin{align*}
a_{n+1}b_{n}+b_{n+1}d_{n}  &  =2b_{n}\cos\frac{2\pi(n+1)}{3}\\
b_{n}c_{n}-a_{n}d_{n}  &  =\frac{b_{n}}{b_{n+1}}\\
\alpha_{n}d_{n}-\beta_{n}b_{n}  &  =\frac{\alpha_{n+1}b_{n}}{b_{n+1}}.
\end{align*}

These equalities imply that%
\[
A_{n}=2\cos\frac{2\pi(n+1)}{3},\quad B_{n}=1,\quad C_{n}=0
\]

In this case, (\ref{dl}) folds into the following:
\begin{equation}
x_{n+1}=2\cos\frac{2\pi n}{3}x_{n}+x_{n-1} \label{c1}%
\end{equation}
which has a coefficient of period 3. The numbers $\alpha_{j},\beta_{j}$ are
readily calculated as
\[
\alpha_{2}=-1,\ \alpha_{3}=2,\ \alpha_{4}=3,\ \beta_{2}=1,\ \beta
_{3}=-1,\ \beta_{4}=-1.
\]

The quadratic equation (\ref{p2qc}) $2r^{2}-4r+1=0$ in this case has two roots
$(2\pm\sqrt{2})/2$. Let $r_{1}=(2-\sqrt{2})/2$ and use (\ref{epr}) to
calculate $r_{2}=1+\sqrt{2}$, $r_{3}=-2+\sqrt{2}$. Since these are units in
$\mathbb{R}$, by Theorem \ref{per} a unitary eigensequence with period 3 is
obtained. The semiconjugate factorization of (\ref{c1}) is readily calculated
and the solution of the factor equation is found to be%
\[
t_{3j+1}=\rho^{j}t_{1},\ t_{3j+2}=-\frac{\rho^{j}t_{1}}{r_{1}},\ t_{3j+3}%
=\frac{\rho^{j}t_{1}}{r_{1}r_{2}},\ j\geq0,\ t_{1}=x_{1}-r_{1}x_{0}%
,\ \rho=-1/(r_{1}r_{2}r_{3}).
\]

Since $\rho=1+\sqrt{2}>1$ it follows that all solutions of (\ref{c1}) with
$t_{1}\not =0$ are unbounded. However, for initial values satisfying
$x_{1}=r_{1}x_{0}$ we have $t_{1}=0$; so $t_{n}=0$ for all $n$. When inserted
in the cofactor equation $x_{n+1}=r_{n+1}x_{n}+t_{n+1}$ this yields%
\begin{equation}
x_{3n}=\frac{(-1)^{n}x_{0}}{\rho^{n}},\ x_{3n+1}=\frac{(-1)^{n}x_{0}r_{1}%
}{\rho^{n}},\ x_{3n+2}=\frac{(-1)^{n}x_{0}r_{1}r_{2}}{\rho^{n}},\ n\geq1.
\label{sps}%
\end{equation}

These special solutions of (\ref{c1}) converge to 0 exponentially for all
$x_{0}$. If $\{x_{n}\}$ is any solution of (\ref{c1}) then the y-component of
the corresponding orbit $\{(x_{n},y_{n})\}$ of the system is given by
(\ref{linyn}). For the special solutions (\ref{sps}), if the coefficient
$a_{n}$ is bounded then $y_{n}\approx-\alpha_{n}/b_{n}$ for large $n$.

\subsection{Linear differential systems and equations}

Many of the results stated above for difference systems and equations have
differential analogs. The essentials are as follows: Consider the system of
linear differential equations%
\begin{equation}
\left\{
\begin{array}
[c]{c}%
x^{\prime}(t)=a(t)x(t)+b(t)y(t)+\alpha(t)\\
y^{\prime}(t)=c(t)x(t)+d(t)y(t)+\beta(t)
\end{array}
\right.  \label{ldfs}%
\end{equation}
with differentiable functions $a(t),b(t),c(t),d(t),\alpha(t),\beta(t)$ defined
on $\mathbb{R}$ or some open subset of it.

If $b(t)\not =0$ for all $t$ then $f(t,x,y)=a(t)x(t)+b(t)y(t)+\alpha(t)$ is a
semi-invertible function that yields%
\begin{equation}
y(t)=\frac{1}{b(t)}[x^{\prime}(t)-a(t)x(t)-\alpha(t)]=h(t,x,x^{\prime}).
\label{eyt}%
\end{equation}

Further, using (\ref{de2}) or by straightforward calculation under suitable
differentiability hypotheses,%
\begin{equation}
x^{\prime\prime}=\left(  a+d+\frac{b^{\prime}}{b}\right)  x^{\prime}+\left(
bc-ad+a^{\prime}-\frac{b^{\prime}a}{b}\right)  x-d\alpha+b\beta-\frac
{b^{\prime}\alpha}{b}+\alpha^{\prime} \label{del}%
\end{equation}
where for brevity we have omitted explicit mention of the variable $t.$ If
$x(t)$ is a solution of (\ref{del}) then the x-component of an orbit of
(\ref{ldfs}) is determined and the y-component is readily found using
(\ref{eyt}).

If $b(t)=0$ for some $t$ but $c(t)\not =0$ for all $t$ then the above
calculations may be repeated with the roles of x and y switched. If both
$b(t)$ and $c(t)$ vanish for some values of $t$ then (\ref{eyt}) is not
applicable and procedures discussed for difference equations can be
implemented here too. These results yield a second-order differential equation
similar to (\ref{del}) that gives the x-component of an orbit of the
differential system (\ref{ldfs}); the y-component is then obtained using the
second equation in (\ref{ldfs}).

If $b(t)=b$ is a constant then $b^{\prime}(t)=0$ and (\ref{del}) reduces to%
\begin{equation}
x^{\prime\prime}=[a(t)+d(t)]x^{\prime}+[bc(t)-a(t)d(t)+a^{\prime
}(t)]x-d(t)\alpha(t)+b\beta(t)+\alpha^{\prime}(t) \label{del0}%
\end{equation}
which is analogous to (\ref{rcom}) in the discrete case. In particular, if
$a(t)=a$ is a constant then the trace and the determinant of the coefficients
matrix of (\ref{ldfs}) can be identified in the coefficients of (\ref{del0}).
The following summarizes the above results for differential systems.

\begin{proposition}
Consider the linear differential system (\ref{ldfs}) in $\mathbb{R}^{2}$.

(a) If $b(t)\not =0$ for all $t$ then the linear differential equation
(\ref{del}) together with the passive equation (\ref{eyt}) constitute a
folding of (\ref{ldfs}).

(b) If $b(t)=b$ is a constant for all $t$ then (\ref{ldfs}) folds into the
linear differential equation (\ref{del0}) which yields the x-component of its
orbits while its y-component is determined either (i) passively from
(\ref{eyt}) if $b\not =0$ or (ii) as a solution of the second equation of
(\ref{ldfs}) if $b=0$.
\end{proposition}

\section{An inverse problem\label{inv}}

Folding a given nonlinear system into a higher order equation does not
necessarily simplify the study of solutions. From a practical point of view, a
significant gain in terms of simplifying the analysis of solutions is
desirable. To address this issue \textit{systematically} we chart a course
backward, from a higher order equation with a desirable property to the system
that yields it through the folding procedure.

In this section we determine and study classes of systems that fold into
difference or differential equations of order 2 with known properties. We
assume that one of the two equations of the system, say, the one specified by
$f,$ is given along with a known function $\phi$ that defines a second-order
equation. Then a function $g$ is determined with the property that the system
with components $f$ and $g$ folds into an equation of order 2 defined by
$\phi.$

This inverse process in particular yields a (non-standard) unfolding of $\phi$
that is based on the given component $f.$ In the special case $f(t,u,v)=v$ we
obtain $g=\phi$ so the corresponding system is just a familiar, standard
unfolding of the equation defined by $\phi.$

\subsection{Difference equations}

Suppose that a function $f$ satisfies condition (\ref{sol}). By (\ref{se2})
the following
\[
f(n+1,w,g(n,u,h(n,u,w)))=\phi(n,u,w)
\]
is a function of $n,u,w.$ Since $f$ is semi-invertible, once again using
(\ref{sol}) we obtain%
\begin{equation}
g(n,u,h(n,u,w))=h(n+1,w,\phi(n,u,w)) \label{idc}%
\end{equation}

Now, suppose that $\phi(n,u,w)$ is prescribed on a set $\mathbb{N\times
}M^{\prime}$ where $M^{\prime}\subset S\times S$ and we seek $g$ that
satisfies (\ref{idc}). Assume that a subset $M$ of $D$ exists with the
property that $f(\mathbb{N\times}M)\times\phi(\mathbb{N\times}M^{\prime
})\subset M^{\prime}.$ For $(n,u,v)\in\mathbb{N\times}M$ define%
\begin{equation}
g(n,u,v)=h(n+1,f(n,u,v),\phi(n,u,f(n,u,v))) \label{idc1}%
\end{equation}

In particular, if $v\in h(\mathbb{N\times}M^{\prime})$ then $g$ above
satisfies (\ref{idc1}). These observations establish the following result.

\begin{theorem}
\label{idn}Let $f$ be a semi-invertible function with $h$ given by
(\ref{sol}). Further, let $\phi$ be a given function on $\mathbb{N\times
}M^{\prime}$. If $g$ is given by (\ref{idc1}) then (\ref{s1}) folds into the
difference equation $s_{n+2}=\phi(n,s_{n},s_{n+1})$ plus a passive equation.
\end{theorem}

As a check, consider $f(u,v)=uv$ and $\phi(u,w)=u(a+bu)$, both independent of
$n$. Then $h(u,w)=w/u$ and (\ref{idc1}) gives%
\[
g(u,v)=\frac{u(a+bu)}{uv}=\frac{a+bu}{v}%
\]

This $g$ yields (\ref{bs1}), as expected. The function $f$ in this example is
separable. In separable cases, explicit expressions are possible with the aid
of (\ref{spsh}). Note that semilinear systems are included in the next result.

\begin{corollary}
Let $(G,\ast)$ be a nontrivial group and $f(n,u,v)=f_{1}(n,u)\ast f_{2}(n,v)$
be separable on $G\times G$ with $f_{2}$ a bijection. If $\phi$ is a given
function on $\mathbb{N\times}G\times G$ and $g$ is given by
\[
g(n,u,v)=f_{2}^{-1}(n+1,\left[  f_{1}(n+1,f_{1}(n,u)\ast f_{2}(n,v))\right]
^{-1}\ast\phi(n,u,f_{1}(n,u)\ast f_{2}(n,v)))
\]
then (\ref{s1}) folds into the difference equation $s_{n+2}=\phi
(n,s_{n},s_{n+1})$ plus a passive equation.
\end{corollary}

The next result yields a class of systems that reduce (effectively) to
first-order difference equations.

\begin{corollary}
\label{cdn}Assume that $f,h$ satisfy the hypotheses of Theorem \ref{idn} and
let $\phi(n,\cdot)$ be a function of one variable for each $n.$ If%
\begin{equation}
g(n,u,v)=h(n+1,f(n,u,v),\phi(n,u)) \label{ief}%
\end{equation}
then (\ref{s1}) folds into the difference equation $s_{n+2}=\phi(n,s_{n})$
whose even terms and odd terms are solutions of the first-order equation%
\[
r_{n+1}=\phi(n,r_{n}).
\]

\end{corollary}

For example, if $f(n,u,v)=(-1)^{n}u+v$ and $\phi(s)$ is a given function
(independent of $n$) then $h(n,u,w)=w-(-1)^{n}u$ and (\ref{ief}) yields%
\[
g(n,u,v)=\phi(u)-(-1)^{n+1}[(-1)^{n}u+v]=\phi(u)+u+(-1)^{n}v
\]

Thus the non-autonomous semilinear system%
\[
\left\{
\begin{array}
[c]{l}%
x_{n+1}=(-1)^{n}x_{n}+y_{n}\\
y_{n+1}=\phi(x_{n})+x_{n}+(-1)^{n}y_{n}%
\end{array}
\right.
\]
folds in the sense of Corollary (\ref{cdn}) into the autonomous, first-order
difference equation $r_{n+1}=\phi(r_{n})$ plus a passive equation. Note that
the above system is the same as (\ref{nas}) with $\psi(u)=\phi(u)+u$.

The next result yields a class of systems that actually reduce to first-order
difference equations.

\begin{corollary}
\label{o1}Assume that $f,h$ satisfy the hypotheses of Theorem \ref{idn} and
let $\phi(n,\cdot)$ be a function of one variable for each $n$. If%
\[
g(n,u,v)=h(n+1,f(n,u,v),\phi(n,f(n,u,v)))
\]
then (\ref{s1}) folds into the difference equation $s_{n+2}=\phi(n,s_{n+1})$
with order 1 plus a passive equation.
\end{corollary}

For example, let $f(u,v)=u/v$ so that $h(u,w)=u/w.$ If $\phi(s)=as+b$ then
define%
\[
g(u,v)=h\left(  f(u,v),af(u,v)+b\right)  =\frac{u/v}{au/v+b}=\frac{u}{au+bv}.
\]

Then the rational system%
\[
\left\{
\begin{array}
[c]{l}%
x_{n+1}=x_{n}/y_{n}\\
y_{n+1}=x_{n}/(ax_{n}+by_{n})
\end{array}
\right.
\]
folds into the first-order linear equation $s_{n+2}=as_{n+1}+b$ plus the
passive equation $y_{n}=s_{n}/s_{n+1}$ with $s_{0}=x_{0}$, $s_{1}=x_{0}%
/y_{0}.$

The next result concerns systems that fold into \textit{autonomous} linear
difference equations of order 2.

\begin{corollary}
\label{lin}Assume that $f,h$ satisfy the hypotheses of Theorem \ref{idn} and
let $\phi(u,w)=\alpha u+\beta w$ be a linear function. If%
\[
g(n,u,v)=h(n+1,f(n,u,v),\alpha u+\beta f(n,u,v))
\]
then (\ref{s1}) folds into the linear autonomous difference equation
$s_{n+2}=\alpha s_{n}+\beta s_{n+1}$ plus a passive equation.
\end{corollary}

The procedure described in this section generates a system that may be
considered an unfolding (non-standard) of the equation $s_{n+2}=\phi
(n,s_{n},s_{n+1}).$ Indeed, if $f(n,u,v)=v$ then $h(n,u,w)=w$ and by
(\ref{idc1}) $g(n,u,v)=\phi(n,u,v),$ i.e., $g=\phi$ so the resulting system is
the standard unfolding (\ref{s2}). In general, for each fixed $\phi$ there
exists a distinct unfolding of $s_{n+2}=\phi(n,s_{n},s_{n+1})$ for every
semi-invertible function $f.$

\subsection{Differential equations}

The inverse problem above has a differential analog. Assume that a function
$f$ is semi-invertible with $h$ given by (\ref{des}). Then the right hand side
of (\ref{de2}) may be written as follows as a function of $t,u,w$%
\begin{equation}
f_{t}(t,u,h(t,u,w))+f(t,u,h(t,u,w))f_{u}(t,u,h(t,u,w))+g(t,u,h(t,u,w))f_{v}%
(t,u,h(t,u,w))=\phi(t,u,w) \label{ic1}%
\end{equation}

Now suppose, inversely, that the function $\phi(t,u,w)$ is prescribed and we
wish to determine a function $g$ such that (\ref{ic1}) is satisfied. Define
\begin{equation}
g(t,u,v)=\frac{1}{f_{v}(t,u,v)}\left[  \phi(t,u,f(t,u,v))-f(t,u,v)f_{u}%
(t,u,v)-f_{t}(t,u,v)\right]  \label{gc}%
\end{equation}
provided that $f_{v}(t,u,v)\not =0.$

Note that unlike (\ref{idc1}), using the chain rule made a second application
of semi-invertibility unnecessary in deriving (\ref{gc}). These observations
establish the following result.

\begin{theorem}
\label{cds}Let $f$ be a semi-invertible function and its semi-inversion $h$ be
given by (\ref{des}) and let $\phi$ be a given function. If $f_{v}%
(t,u,v)\not =0$ for all $t,u,v$ and $g$ is given by (\ref{gc}) then system
(\ref{ds1}) folds into the following second-order ordinary differential
equation (plus a passive equation)%
\[
s^{\prime\prime}=\phi(t,s,s^{\prime})
\]

\end{theorem}

To illustrate, consider $f(u,v)=tu^{2}-v$ as in (\ref{des1}) and let
$\phi(u,w)=\alpha u.$ Then $f_{t}(u,v)=u^{2}$, $f_{u}(u,v)=2tu$ and
$f_{v}(u,v)=-1$ so from (\ref{gc})%
\[
g(u,v)=-[\alpha u-2tu(tu^{2}-v)-u^{2}]=-\alpha u+u^{2}+2t^{2}u^{3}-2tuv
\]

If $a=-\alpha$ then $g(u,v)=au+u^{2}+2t^{2}u^{3}-2tuv$ and system (\ref{des1})
is obtained.

Theorem \ref{cds} can be used to categorize classes of differential systems
that fold into second-order, nonlinear differential equations with known
properties. For example, we may determine systems that fold into the equations
of Lienard, van der Pol or Duffing; see, e.g., \cite{Draz} or \cite{Strg}. In
particular, the periodically forced, special case of Duffing's equation%
\[
x^{\prime\prime}+bx^{\prime}+kx^{3}=A\sin\omega t
\]
which generates chaotic flows has a different (non-standard) unfolding for
every choice of $f(t,u,v)$ with $f_{v}\not =0.$ For instance, if
$f(t,u,v)=\alpha u+\beta v$ for constants $\alpha,\beta$ with $\beta\not =0$
then with $\phi(t,u,w)=A\sin\omega t-ku^{3}-bw,$ (\ref{gc}) yields%
\begin{align*}
g(t,u,v)  &  =\frac{1}{\beta}\left[  A\sin\omega t-ku^{3}-b(\alpha u+\beta
v)-\alpha(\alpha u+\beta v)\right] \\
&  =\frac{A}{\beta}\sin\omega t-\frac{k}{\beta}u^{3}-\frac{\alpha}{\beta
}(\alpha+b)u-(\alpha+b)v
\end{align*}

In particular, if we set $\alpha=-b$ and $\beta=1$ then the system%
\[
\left\{
\begin{array}
[c]{l}%
x^{\prime}=-bx+y\\
y^{\prime}=-kx^{3}+A\sin\omega t
\end{array}
\right.
\]
folds into Duffing's equation plus the passive equation $y=x^{\prime}+bx;$
other variations are clearly possible.

\begin{remark}
Analogs of Corollaries (\ref{cdn}), (\ref{o1}) and (\ref{lin}) for
differential systems are also true but we do not list them explicitly.
\end{remark}

\section{Folding larger systems into equations}

Folding a system of $k$ equations into an equation of order $k\geq3$ using the
method of the preceding sections is also possible if suitable inversions
exist. The process of folding a system of $k$ equations may in principle be
viewed as iterating what was discussed above. The larger the value of $k$, the
longer and more complex the $k$-th order equation needs to be in order to
include all the information that is contained in the system.

Linear systems illustrate this situation: while a system of 2 linear equations
is defined by 4 parameters (the coefficient matrix has 4 entries) a system of
3 linear equations requires 9 parameters. The increased amount of information,
both algebraically (dealing with a $3\times3$ matrix) and analytically (orbits
exist in 3 dimensions rather than 2) must be configured within a single
equation. For \textit{nonlinear} equations potential technical difficulties
related to inversions need to be considered.

In certain cases where the higher order equation that is derived from the
system is much simpler than the system itself, folding yields a practical
advantage. We discuss examples of such systems in this section, including a
special type of difference system with $k$ equations that readily folds into
an equation of order $k$ without requiring any type of inversion.

\subsection{Folding systems of 3 equations}

In this section, we discuss folding systems of 3 difference or differential
equations into higher order equations. Let $S$ be a nonempty set and let
$f_{1},f_{2},f_{3}:\mathbb{N}\times D\rightarrow S$ be given functions where
$D\subset S\times S\times S$. Consider the system%
\begin{equation}
\left\{
\begin{array}
[c]{l}%
x_{1,n+1}=f_{1}(n,x_{1,n},x_{2,n},x_{3,n})\\
x_{2,n+1}=f_{2}(n,x_{1,n},x_{2,n},x_{3,n})\\
x_{3,n+1}=f_{3}(n,x_{1,n},x_{2,n},x_{3,n})
\end{array}
\right.  \label{3eq}%
\end{equation}

Starting from an initial point $(x_{1,0},x_{2,0},x_{3,0})\in D,$ iteration
generates points $(x_{1,n},x_{2,n},x_{3,n})$ of an orbit of (\ref{3eq}) for as
long as the orbit stays within the domain $D$ of the system. Assume that for
every $n$ the following \textit{partial inversion} holds%
\begin{equation}
r=f_{1}(n,u,v,w)\Rightarrow w=\phi_{1}(n,u,v,r) \label{zsmi}%
\end{equation}

Then the first equation of (\ref{3eq}) yields%
\[
x_{3,n}=\phi_{1}(n,x_{1,n},x_{2,n},x_{1,n+1}),\quad n=1,2,\ldots
\]

Now%
\begin{align*}
x_{1,n+2}  &  =f_{1}(n+1,x_{1,n+1},x_{2,n+1},x_{3,n+1})\\
&  =f_{1}(n+1,x_{1,n+1},f_{2}(n,x_{1,n},x_{2,n},x_{3,n}),f_{3}(n,x_{1,n}%
,x_{2,n},x_{3,n}))\\
&  =f_{1}(n+1,x_{1,n+1},f_{2}(n,x_{1,n},x_{2,n},\phi_{1}(n,x_{1,n}%
,x_{2,n},x_{1,n+1})),f_{3}(n,x_{1,n},x_{2,n},\phi_{1}(n,x_{1,n},x_{2,n}%
,x_{1,n+1})))
\end{align*}

The last expression does not involve the variable $x_{3,n}$ so we use it to
define a new function%
\[
f_{1}^{(1)}(n,u,v,r)=f_{1}(n+1,r,f_{2}(n,u,v,\phi_{1}(n,u,v,r)),f_{3}%
(n,u,v,\phi_{1}(n,u,v,r)))
\]

Thus, $x_{1,n+2}=f_{1}^{(1)}(n,x_{1,n},x_{2,n},x_{1,n+1}).$ Next, assume that
$\bar{f}$ has the following partial inversion%
\begin{equation}
s=f_{1}^{(1)}(n,u,v,r)\Rightarrow v=\phi_{2}(n,u,r,s) \label{fbr}%
\end{equation}

Then $x_{2,n}=\phi_{2}(n,x_{1,n},x_{1,n+1},x_{1,n+2})$ and the following is
obtained:%
\begin{equation}
x_{1,n+3}=f_{1}^{(1)}(n+1,x_{1,n+1},x_{2,n+1},x_{1,n+2})=f_{1}^{(2)}%
(n,x_{1,n},x_{1,n+1},x_{1,n+2}), \label{sde3}%
\end{equation}

where%
\[
f_{1}^{(2)}(n,u,r,s)=f_{1}^{(1)}(n+1,r,f_{2}(n,u,\phi_{2}(n,u,r,s),\phi
_{1}(n,u,\phi_{2}(n,u,r,s),r)),s).
\]

If $\{x_{1,n}\}$ is a solution of (\ref{sde3}) with initial values%
\begin{align*}
x_{1,0},\ x_{1,1}  &  =f_{1}(0,x_{1,0},x_{2,0},x_{3,0}),\ x_{1,2}%
=f_{1}(1,x_{1,1},x_{2,1},x_{3,1})\\
\text{where\quad}x_{2,1}  &  =f_{2}(0,x_{1,0},x_{2,0},x_{3,0}),\ x_{3,1}%
=f_{3}(0,x_{1,0},x_{2,0},x_{3,0})
\end{align*}
then the second and third components of an orbit $\{(x_{1,n},x_{2,n}%
,x_{3,n})\}$ are obtained passively via the partial inversion equations:%
\begin{equation}
x_{2,n}=\phi_{2}(n,x_{1,n},x_{1,n+1},x_{1,n+2}),\quad x_{3,n}=\phi
_{1}(n,x_{1,n},x_{2,n},x_{1,n+1}) \label{3p}%
\end{equation}
without needing to solve additional difference equations.

As in Definition \ref{fld}, we say that the third-order equation (\ref{sde3})
together with the two passive equations (\ref{3p}) constitute a folding of
(\ref{3eq}).

To illustrate the above process, we apply it to the system%
\begin{equation}
\left\{
\begin{array}
[c]{l}%
x_{1,n+1}=ax_{2,n}+x_{3,n}\\
x_{2,n+1}=bx_{1,n}+cx_{3,n}\\
x_{3,n+1}=\rho(x_{1,n})+dx_{2,n}%
\end{array}
\right.  \label{i3es}%
\end{equation}
where $a,b,c,d\in\mathbb{R}$ and $\rho:\mathbb{R\rightarrow R}$ is a given
function. From the first equation of (\ref{i3es}),%
\begin{equation}
x_{3,n}=x_{1,n+1}-ax_{2,n} \label{x3n}%
\end{equation}

This is the function $\phi_{1}$ in the above algorithm. Next, shifting the
index by 1 in the first equation of (\ref{i3es}) and using substitutions%
\begin{equation}
x_{1,n+2}=ax_{2,n+1}+x_{3,n+1}=abx_{1,n}+acx_{3,n}+x_{3,n+1}=acx_{1,n+1}%
+abx_{1,n}+\rho(x_{1,n})+(d-a^{2}c)x_{2,n} \label{x1n2}%
\end{equation}

This equation is easily solved for $x_{2,n}$ to obtain the function $\phi_{2}$
in the above algorithm. Next, shift the index once more by 1 and use
substitutions to obtain%
\begin{align}
x_{1,n+3}  &  =ax_{2,n+2}+x_{3,n+2}=abx_{1,n+1}+acx_{3,n+1}+\rho
(x_{1,n+1})+dx_{2,n+1}\nonumber\\
&  =abx_{1,n+1}+\rho(x_{1,n+1})+ac\rho(x_{1,n})+acdx_{2,n}+d(bx_{1,n}%
+cx_{3,n})\nonumber\\
&  =abx_{1,n+1}+\rho(x_{1,n+1})+ac\rho(x_{1,n})+acdx_{2,n}+bdx_{1,n}%
+cd(x_{n,n+1}-ax_{2,n})\nonumber\\
&  =(ab+cd)x_{1,n+1}+\rho(x_{1,n+1})+bdx_{1,n}+ac\rho(x_{1,n}) \label{i3de}%
\end{align}

Since the terms containing $x_{2,n}$ cancel out in the substitution process,
it is not necessary to substitute $\phi_{2}$ for $x_{2,n}$ in this case to
obtain (\ref{i3de}). But deriving $\phi_{2}$ is still necessary to calculate
the orbits of (\ref{i3es}) once a solution $\{x_{1,n}\}$ of (\ref{i3de}) is
determined. We distinguish two cases: If $d\not =a^{2}c$ then the passive
equation%
\[
x_{2,n}=\frac{1}{d-a^{2}c}(x_{1,n+2}-acx_{1,n+1}-abx_{1,n}-\rho(x_{1,n}))
\]
yields the component $x_{2,n}$ and then (\ref{x3n}) gives the third component
$x_{3,n}.$ If $d=a^{2}c$ then the second-order difference equation
(\ref{x1n2}) reduces to%
\begin{equation}
x_{1,n+2}=acx_{1,n+1}+abx_{1,n}+\rho(x_{1,n}) \label{o2n3}%
\end{equation}
indicating that (\ref{i3es}) is not strictly 3-dimensional in this case. The
second component $x_{2,n}$ may be found from the system by eliminating
$x_{3,n}$ from the first two equations in (\ref{i3es}) to obtain the
difference equation%
\[
x_{2,n+1}=bx_{1,n}+cx_{3,n}=-acx_{2,n}+bx_{1,n}+cx_{1,n+1}%
\]

With the sequence $\{x_{1,n}\}$ being given as a solution of (\ref{o2n3}), the
above equation is a linear, nonhomogeneous first-order equation that yields
the component $x_{2,n}$ as a solution. Then (\ref{x3n}) yields $x_{3,n}$.
Alternatively, we may find $x_{3,n}$ first by eliminating $x_{2,n}$ from the
first and the third equations in (\ref{i3es}) to obtain a linear difference
equation that yields $x_{3,n}$ as a solution. Then $x_{2,n}$ maybe found using
any equation of the system.

In principle, the algorithm described above can be extended to systems
containing any number $k$ of equations; we discuss a general folding algorithm
for difference systems in the next section.

Systems of 3 differential equations also fold into ordinary differential
equations using the analog of the above algorithm in which higher order
derivatives replace index shifts. We do not enter the details for the
differential case here but illustrate the process using the following analog
of system (\ref{i3es})
\begin{equation}
\left\{
\begin{array}
[c]{l}%
x_{1}^{\prime}(t)=ax_{2}(t)+x_{3}(t)\\
x_{2}^{\prime}(t)=bx_{1}(t)+cx_{3}(t)\\
x_{3}^{\prime}(t)=\rho(x_{1}(t))+dx_{2}(t)
\end{array}
\right.  \label{des3e}%
\end{equation}

The first equation of (\ref{des3e}) yields%
\[
x_{3}(t)=x_{1}^{\prime}(t)-ax_{2}(t).
\]

Now taking the derivative of the first equation of (\ref{des3e}) and using
substitutions%
\begin{align*}
x_{1}^{\prime\prime}(t)  &  =ax_{2}^{\prime}(t)+x_{3}^{\prime}(t)=abx_{1}%
(t)+acx_{3}(t)+\rho(x_{1}(t))+dx_{2}(t)\\
&  =abx_{1}(t)+\rho(x_{1}(t))+acx_{1}^{\prime}(t)+(d-a^{2}c)x_{2}(t)
\end{align*}

The above equation may be solved for $x_{2}(t)$ if $d\not =a^{2}c$. Next, take
derivatives once more and use substitutions to obtain%
\begin{align*}
x_{1}^{\prime\prime\prime}(t)  &  =ax_{2}^{\prime\prime}(t)+x_{3}%
^{\prime\prime}(t)=abx_{1}^{\prime}(t)+acx_{3}^{\prime}(t)+\rho^{\prime}%
(x_{1}(t))x_{1}^{\prime}(t)+dx_{2}^{\prime}(t)\\
&  =abx_{1}^{\prime}(t)+\rho^{\prime}(x_{1}(t))x_{1}^{\prime}(t)+ac\rho
(x_{1}(t))+acdx_{2}(t)+dbx_{1}(t)+dcx_{3}(t)\\
&  =abx_{1}^{\prime}(t)+\rho^{\prime}(x_{1}(t))x_{1}^{\prime}(t)+ac\rho
(x_{1}(t))+acdx_{2}(t)+bdx_{1}(t)+dcx_{1}^{\prime}(t)-dcax_{2}(t)\\
&  =(\rho^{\prime}(x_{1}(t))+ab+cd)x_{1}^{\prime}(t)+ac\rho(x_{1}%
(t))+bdx_{1}(t)
\end{align*}

The last ordinary differential equation above is what system (\ref{des3e})
folds into. If $x_{1}(t)$ is a solution of this third-order differential
equation then $x_{2}(t)$ and $x_{3}(t)$ are determined in a passive way that
is analogous to the determination of $x_{2,n}$ and $x_{3,n}$ for the
difference system (\ref{i3es}).

Of interest is a comparison (in the autonomous case) of the above algorithmic
or iterative method with the method that is discussed in \cite{Eich}. The
latter method applies to a general class of differential systems of three
equations and yields a jerk function provided certain conditions hold.
Jacobian inversion (see \cite{Abhy}) is used (possibly over-used) to derive a
type of system that is amenable to being folded in such a way that the flows
of the system are represented by the scalar equation. A system of the
following type is derived%
\begin{equation}
\left\{
\begin{array}
[c]{l}%
x_{1}^{\prime}(t)=\eta_{1}(x_{1}(t))+b_{11}x_{2}(t)+b_{13}x_{3}(t)+c_{1}\\
x_{2}^{\prime}(t)=\eta_{2}(x_{1}(t),x_{2}(t),x_{3}(t))+b_{21}x_{1}%
(t)+b_{22}x_{2}(t)+b_{23}x_{3}(t)+c_{2}\\
x_{3}^{\prime}(t)=\eta_{3}(x_{1}(t),x_{2}(t),x_{3}(t))+b_{31}x_{1}%
(t)+b_{32}x_{2}(t)+b_{33}x_{3}(t)+c_{3}%
\end{array}
\right.  \label{eich}%
\end{equation}
where constants $b_{ij}$, $c_{i}$ and functions $\eta_{i}$ are arbitrary
except for the following conditions%
\begin{align}
b_{12}\eta_{2}(x_{1}(t),x_{2}(t),x_{3}(t))+b_{13}\eta_{2}(x_{1}(t),x_{2}%
(t),x_{3}(t))  &  =\xi(x_{1}(t),b_{12}x_{2}(t)+b_{13}x_{3}(t))\label{eich1}\\
b_{12}^{2}b_{23}-b_{13}^{2}b_{32}+b_{12}b_{13}(b_{33}-b_{22})  &  \not =0
\label{eich2}%
\end{align}
with $\xi$ an arbitrary function. Condition (\ref{eich2}) is needed to ensure
that (\ref{eich}) folds to a third-order (jerk) equation and not to a
second-order one; see \cite{Eich}. Condition (\ref{eich1}) ensures the
elimination of variables $x_{2}$ and $x_{3}$ in a relatively general setting.

The linear nature of the first equation in the variables $x_{2}$ and $x_{3}$
in (\ref{eich}) is due to the derivation in \cite{Eich}; it makes partial
inversion and solving globally for $x_{3}$ (or $x_{2}$) easy. However, this
restriction is not strictly necessary; for example, a separable form with
invertible functions would also work in the algorithmic case as in the
following system%
\begin{equation}
\left\{
\begin{array}
[c]{l}%
x_{1}^{\prime}(t)=\eta(x_{1}(t))+\zeta(x_{3}(t))\\
x_{2}^{\prime}(t)=g(x_{1}(t),x_{2}(t),x_{3}(t))\\
x_{3}^{\prime}(t)=\rho(x_{1}(t))+ax_{2}(t)
\end{array}
\right.  \label{fa}%
\end{equation}
where $a\not =0$, $\eta,g$, $\rho$ are arbitrary functions adn $\zeta$ is a
bijection of $\mathbb{R}$. If these maps are nonlinear then the above system
is not of type (\ref{eich}) but may be folded using the folding algorithm. The
first equation of (\ref{fa}) yields $x_{3}(t)=\zeta^{-1}(x_{1}^{\prime
}(t)-\eta(x_{1}(t)))$ so
\[
x_{1}^{\prime\prime}(t)=\eta^{\prime}(x_{1}(t))x_{1}^{\prime}(t)+\zeta
^{\prime}(x_{3}(t))x_{3}^{\prime}(t)=\eta^{\prime}(x_{1}(t))x_{1}^{\prime
}(t)+\zeta^{\prime}(\zeta^{-1}(x_{1}^{\prime}(t)-\eta(x_{1}(t))))[\rho
(x_{1}(t))+ax_{2}(t)]
\]

The above equation is readily solved for $x_{2}(t)$ in terms of the variables
$x_{1}(t),x_{1}^{\prime}(t)$ and $x_{1}^{\prime\prime}(t)$ so $x_{1}%
^{\prime\prime\prime}(t)$ can then be also expressed in terms of
$x_{1}(t),x_{1}^{\prime}(t)$ and $x_{1}^{\prime\prime}(t).$ These observations
suggest that the algorithmic approach is more general, if not the most succinct.

\subsection{A folding algorithm for difference systems}

In this section we extend the folding process for a system of 3 equations
discussed the preceding section to a general folding algorithm for systems of
$k$ difference equations. The process is similar, in essence, for differential
systems but we leave the differential case for future considerations. In each
step of the iterative process a \textquotedblleft space" variable $x_{j,n}$ is
converted to a \textquotedblleft temporal extension" of $x_{1,n}$. Iteration
stops when the $k-1$ space variables $x_{2,n},\ldots,x_{k,n}$ have been
converted to the temporal extensions $x_{1,n+1},\ldots,x_{1,n+k}$ (or possibly
earlier). This algorithm arbitrarily centers on the variable $x_{1,n}$ for
convenience, although in principle any one of $x_{2,n},\ldots,x_{k,n}$ can be
used instead if it simplifies calculations.

We begin with the recursive system%
\begin{equation}
x_{i,n+1}=f_{i}(n,x_{1,n},x_{2,n},\ldots,x_{k,n}),\quad i=1,2,\ldots,k
\label{sysk}%
\end{equation}
of $k$ first-order difference equations in $k$ variables $x_{1,n}%
,x_{2,n},\ldots,x_{k,n}$.

\begin{enumerate}
\item \textit{Shift} the indices of the first equation ($i=1$) in (\ref{sysk})
by 1 to obtain%
\[
x_{1,n+2}=f_{1}(n+1,x_{1,n+1},x_{2,n+1},\ldots,x_{k,n+1})
\]

\begin{enumerate}
\item Use (\ref{sysk}) to \textit{substitute}
\[
x_{j,n+1}=f_{j}(n,x_{1,n},x_{2,n},\ldots,x_{k,n}),\quad j=2,\ldots,k.
\]

\item Assume that $f_{1}$ has the following \textit{partial inversion}%
\[
s_{1}=f_{1}(n,u_{1},u_{2},\ldots,u_{k})\Rightarrow u_{k}=\phi_{1}%
(n,s_{1},u_{1},u_{2},\ldots,u_{k-1})
\]

\item \textit{Substitute} $\phi_{1}(n,s_{1},u_{1},\ldots,u_{k-1})$ for $u_{k}$
in $f_{j}(n,u_{1},u_{2},\ldots,u_{k})$ for $j\geq2$ to define%
\begin{gather*}
f_{1}^{(1)}(n,s_{1},u_{1},\ldots,u_{k-1})=f_{1}(n+1,s_{1},u_{2}^{(1)}%
,\ldots,u_{k}^{(1)})\\
\text{where}\text{:\quad\ }u_{j}^{(1)}=f_{j}(n,u_{1},\ldots,u_{k-1},\phi
_{1}(n,s_{1},u_{1},\ldots,u_{k-1}))
\end{gather*}

Since $u_{j}^{(1)}$ stands for $x_{j,n+1}$ for $j\geq2$ we may write%
\begin{equation}
x_{1,n+2}=f_{1}^{(1)}(n,x_{1,n+1},x_{2,n},\ldots,x_{k-1,n}). \label{sysk1}%
\end{equation}

Notice that the variable $x_{k,n}$ has been replaced with the new variable
$x_{1,n+1}$.
\end{enumerate}

\item \textit{Shift} the indices in (\ref{sysk1}) by 1 to obtain%
\[
x_{1,n+3}=f_{1}^{(1)}(n+1,x_{1,n+2},x_{2,n+1},\ldots,x_{k-1,n+1})
\]

\begin{enumerate}
\item Use (\ref{sysk}) to \textit{substitute} $f_{j}(n,x_{1,n},x_{2,n}%
,\ldots,x_{k,n})$ for $x_{j,n+1}$ for $j=2,\ldots,k-1,$ as in Step 1a above.

\item Assume that $f_{1}^{(1)}$ has the following \textit{partial inversion}%
\[
s_{2}=f_{1}^{(1)}(n,s_{1},u_{1},\ldots,u_{k-1})\Rightarrow u_{k-1}=\phi
_{2}(n,s_{2},s_{1},u_{1},\ldots,u_{k-2})
\]

\item \textit{Substitute} $\phi_{1}(n,s_{1},u_{1},\ldots,u_{k-1})$ for $u_{k}$
in $f_{j}(n,u_{1},\ldots,u_{k-1},u_{k})$ and then \textit{substitute} the
quantity $\phi_{2}(n,s_{2},s_{1},u_{1},\ldots,u_{k-2})$ for $u_{k-1}$ in both
of the preceding quantities to obtain%
\[
f_{1}^{(2)}(n,s_{2},s_{1},u_{1},\ldots,u_{k-2})=f_{1}^{(1)}(n+1,s_{2}%
,u_{2}^{(2)},\ldots,u_{k-1}^{(2)})
\]

where, abbreviating $\phi_{2}(n,s_{2},s_{1},u_{1},\ldots,u_{k-2})$ by
$\phi_{2}$ to save space, we have
\[
u_{j}^{(2)}=f_{j}(n,u_{1},\ldots,u_{k-2},\phi_{2},\phi_{1}(n,s_{1}%
,u_{1},\ldots,u_{k-2},\phi_{2}))\quad j=2,\ldots,k-1.
\]

Since $u_{j}^{(2)}$ stands for $x_{j,n+2}$ for $j\geq2$ we may write%
\[
x_{1,n+3}=f_{1}^{(2)}(n,x_{1,n+2},x_{1,n+1},x_{2,n},\ldots,x_{k-2,n}).
\]

In this step, the variables $x_{k,n}$ and $x_{k-1,n}$ have been replaced with
the new variables $x_{1,n+1},x_{1,n+2}$.
\end{enumerate}
\end{enumerate}

The above procedure repeats as long as inversions are possible. The $i$-th
step of the algorithm is as follows:

\begin{quotation}
\textbf{i.} With $x_{1,n+i}=f_{1}^{(i-1)}(n,x_{1,n+i-1},x_{2,n},\ldots
,x_{k-i+1,n})$ obtained in Step $i-1$, shift indices by 1 to get%
\[
x_{1,n+i+1}=f_{1}^{(i-1)}(n+1,x_{1,n+i},x_{2,n+1},\ldots,x_{k-i,n+1})
\]

\quad(\textbf{a}) Use (\ref{sysk}) to \textit{substitute} $f_{j}%
(n,x_{1,n},x_{2,n},\ldots,x_{k,n})$ for $x_{j,n+1}$ for $j=2,\ldots,k-1,$ as
in Steps 1a, 2a above.

\quad(\textbf{b}) Assume that $f_{1}^{(i-1)}$ has the following
\textit{partial inversion}%
\[
s_{i}=f_{1}^{(i-1)}(n,s_{i-1},\ldots,s_{1},u_{1},\ldots,u_{k-i+1})\Rightarrow
u_{k-i+1}=\phi_{i}(n,s_{i},\ldots,s_{1},u_{1},\ldots,u_{k-i})
\]

\quad(\textbf{c}) A series of \textit{substitutions} now eliminates the
variables $u_{k+j}$ for $j=0,1,\ldots,k-i+1$ as in 2c. Substitute $\phi
_{1}(n,s_{1},u_{1},\ldots,u_{k-1})$ for $u_{k}$ in $f_{j}(n,u_{1}%
,\ldots,u_{k-1},u_{k})$ to eliminate $u_{k}.$ Then substitute $\phi
_{2}(n,s_{2},s_{1},u_{1},\ldots,u_{k-2})$ for $u_{k-1}$ in the expression just
obtained to remove $u_{k},u_{k-1}.$ Continuing, we obtain a function
$f_{1}^{(i)}(n,s_{i},\ldots,s_{1},u_{1},\ldots,u_{k-i})$ \ that is independent
of $u_{k},u_{k-1},\ldots,u_{k-i+1}.$
\end{quotation}

\begin{remark}
Execution of the above algorithm may stop before $k-1$ steps if $f_{1}^{(i)}$
is independent of $u_{2},\ldots,u_{k-i}$ for some $i<k-1$, i.e.,
\[
f_{1}^{(i)}(n,s_{i},\ldots,s_{1},u_{1},\ldots,u_{k-i})=f_{1}^{(i)}%
(n,s_{i},\ldots,s_{1},u_{1}).
\]
An extreme example is the following system%
\begin{align}
x_{1,n+1}  &  =f_{1}(n,x_{1,n},x_{2,n},\ldots,x_{k,n}),\label{1cp}\\
x_{i,n+1}  &  =f_{i}(n,x_{1,n}),\quad i=2,\ldots,k\nonumber
\end{align}
where%
\[
x_{1,n+2}=f_{1}(n+1,x_{1,n+1},x_{2,n+1},\ldots,x_{k,n+1})=f_{1}(n+1,x_{1,n+1}%
,f_{2}(n,x_{1,n}),\ldots,f_{k}(n,x_{1,n}))
\]
is already a difference equation that does not involve any of the variables
$x_{2,n},\ldots,x_{k,n}$ (this system is a special case of the system
discussed in the next section).
\end{remark}

In general, partial inversions may be needed in some steps of the above
algorithm and not in other steps. The folding algorithm above may also stop if
$f_{1}^{(i)}$ has no known partial inversions for some $i<k.$ In this case
using a different variable $x_{j,k}$ for $j>1$ may result in a successful
execution of the folding algorithm.

The preceding discussion suggests a relevant notion that may be associated
with systems to gauge the degree of interdependence of variables.

\begin{definition}
Assume that a system is foldable into a higher order difference (or
differential) equation. Then the order $\kappa$ of the higher order equation
is is the \textbf{interdependence degree (i.d.)} of the system. If $\kappa=k$
then the foldable system is\textbf{ fully interdependent} (or fully coupled).
\end{definition}

In the above language, a completely controllable system in control theory is
fully interdependent. We note that $1\leq\kappa\leq k$ for a foldable system.
The i.d. of a system can often be determined through the use of the folding
algorithm even when explicit formulas for the folding equation (or some of the
partial inversions needed to calculate it) are not known or do not exist;
existence proofs may suffice to carry the algorithm forward and non-existence
arguments may be used to terminate it before $k$ steps.

For a first-order equation $x_{n+1}=f(n,x_{n})$ we define $\kappa=1$ and for
an \textit{uncoupled system }of first-order equations
\[
x_{i,n+1}=f_{i}(n,x_{i,n}),\quad i=1,\ldots,k
\]
with $k\geq2$ we define $\kappa=0.$ Both of these are examples of
\textit{unfoldable} systems since the folding process cannot initiate but
their i.d. values are intuitively obvious.

For some foldable systems, $1\leq\kappa<k;$ e.g., for system (\ref{1cp})
$\kappa=2$ but $k$ may be any integer greater than 1. In general, the value of
$\kappa$ depends not only on the component functions $f_{i}$ but also on the
values of the parameters involved. For example, for (\ref{i3es}) $\kappa=2$ if
$d=a^{2}c$ but $\kappa=3$ otherwise. Such a system has a \textit{variable
i.d.} depending on the parameters. The same conclusions hold for the
differential system (\ref{des3e}). Similarly, for a linear system of two
equations, $\kappa=1$ if the coefficients matrix has determinant 0 and
$\kappa=2$ otherwise; see, e.g., (\ref{rcom}). On the other hand, a linear
system whose coefficients matrix is diagonal with nonzero entries has nonzero
determinant and yet, it is an uncoupled system with $\kappa=0$.

Finally, if a particular system splinters into two or more blocks or
subsystems with disjoint sets of variables for a range of parameter values
then such a system has an undefined i.d. for certain parameter ranges. Of
course, each of the subsystems containing more than one variable may be
separately folded in such cases.

\subsection{Folding without inversions}

There are systems that fold into equations without requiring any partial
inversions. This feature reduces the difficulty of calculations considerably
in nonlinear systems. Consider the following difference system%
\begin{equation}
\left\{
\begin{array}
[c]{l}%
x_{1,n+1}=f_{1}(n,x_{1,n},x_{2,n},x_{3,n},\ldots,x_{k,n})\\
x_{2,n+1}=f_{2}(n,x_{1,n},x_{3,n},\ldots,x_{k,n})\\
\vdots\\
x_{k-1,n+1}=f_{k-1}(n,x_{1,n},x_{k,n})\\
x_{k,n+1}=f_{k}(n,x_{1,n})
\end{array}
\right.  \label{ske}%
\end{equation}
or more concisely,%
\begin{align*}
x_{i,n+1}  &  =f_{i}(n,x_{1,n},x_{i+1,n},x_{i+2,n},\ldots,x_{k,n}),\quad
i=1,2,\ldots,k-1\\
x_{k,n+1}  &  =f_{k}(n,x_{1,n}).
\end{align*}

Here we\textit{ do not assume }that the functions $f_{i}$ have partial
inversions that yield any of the variables\textit{ }$x_{1,n},x_{2,n}%
,\ldots,x_{k,n}$. In this case, rather than using the folding algorithm of the
last section, it is more expedient to proceed as follows. Shift the indices of
the first equation in (\ref{ske}) by $k-1$ to obtain
\begin{equation}
x_{1,n+k}=f_{1}(n+k-1,x_{1,n+k-1},x_{2,n+k-1},x_{3,n+k-1},\ldots,x_{k,n+k-1})
\label{spde}%
\end{equation}

Each component in $f_{1}$ above is determined from the system itself as%
\begin{align}
x_{i,n+k-1}  &  =f_{i}(n+k-2,x_{1,n+k-2},x_{i+1,n+k-2},\ldots,x_{k,n+k-2}%
),\quad i=1,2,\ldots,k-1\nonumber\\
x_{k,n+k-1}  &  =f_{k}(n+k-2,x_{1,n+k-2}) \label{sdk}%
\end{align}

This step gives the last component $x_{k,n+k-1}$ in (\ref{spde}) exclusively
in terms of $x_{1,n+k-2}.$ The remaining components $x_{i,n+k-1}$ for
$i=2,3,\ldots k-1$ are found recursively in terms of $x_{1,n+k-j}$ for
suitable values of $j.$ If $i=k-1$ then%
\begin{align}
x_{k-1,n+k-1}  &  =f_{k-1}(n+k-2,x_{1,n+k-2},x_{k,n+k-2})\nonumber\\
&  =f_{k-1}(n+k-2,x_{1,n+k-2},f_{k}(n+k-3,x_{1,n+k-3})) \label{sdk1}%
\end{align}

Although this expression is more complicated, it is built up from (\ref{sdk})
since the term $x_{k,n+k-2}$ is essentially that in (\ref{sdk}) but with index
reduced by 1. Similarly,%
\begin{align*}
x_{k-2,n+k-1}  &  =f_{k-2}(n+k-2,x_{1,n+k-2},x_{k-1,n+k-2},x_{k,n+k-2})\\
&  =f_{k-2}(n+k-2,x_{1,n+k-2},f_{k-1}(n+k-3,x_{1,n+k-3},x_{k,n+k-3}%
),f_{k}(n+k-3,x_{1,n+k-3}))\\
&  =f_{k-2}(n+k-2,x_{1,n+k-2},f_{k-1}(n+k-3,x_{1,n+k-3},f_{k}%
(n+k-4,x_{1,n+k-4})),\\
&  \hspace{4in}f_{k}(n+k-3,x_{1,n+k-3}))
\end{align*}
in which the terms $x_{k-1,n+k-2}$ and $x_{k,n+k-2}$ are again given by
(\ref{sdk1}) and (\ref{sdk}), respectively but with indices reduced by 1 in
each case. In this way, we determine all of the components in (\ref{spde})
recursively by going backward. The expressions become more involved as the
value of $i$ decreases to 2 but the algorithm is well-defined.

Let $\{x_{1,n}\}$ be a solution of (\ref{spde}) in which the terms
$x_{2,n+k-1},x_{3,n+k-1},\ldots,x_{k,n+k-1}$ are replaced by appropriate
expressions in terms of $x_{1,n+k-j}$ as outlined above so that (\ref{spde})
is a difference equation of order $k.$ Then we may calculate $x_{i,n}$ for
$i=2,\ldots,k$ using the system. Given the initial values $x_{i,0}$ for
$i=1,\ldots,k$ we first calculate the initial values $x_{1,\,j}$,
$j=1,\ldots,k$ for (\ref{spde}) as follows:%
\[
x_{1,\,j}=f_{1}(j-1,x_{1,j-1},x_{2,j-1},\ldots,x_{k,j-1}),\quad x_{i,\,j-1}%
=f_{i}(j-2,x_{1,j-2},x_{i+1,j-2},\ldots,x_{k,j-2})
\]
for $j=2,\ldots,k$\ and of course, $x_{1,\,1}=f_{1}(0,x_{1,0},x_{2,0}%
,\ldots,x_{k,0}).$ Next, once $\{x_{1,n}\}$ is determined for all $n$ from
(\ref{spde}) we calculate $x_{i,n}$ for $i=2,\ldots,k$ in a backward fashion
as follows:%
\begin{align*}
x_{k,n}  &  =f_{k}(n-1,x_{1,n-1}),\quad n\geq1\\
x_{k-1,n}  &  =f_{k-1}(n-1,x_{1,n-1},x_{k,n-1}),\\
&  \vdots\\
x_{2,n}  &  =f_{2}(n-1,x_{1,n-1},x_{3,n-1},\ldots,x_{k,n-1})
\end{align*}

Thus, an orbit of (\ref{ske}) may be determined from a solution of
(\ref{spde}) without having to solve additional difference equations or to
find partial inversions.

To illustrate the above procedure, consider the autonomous system%
\begin{equation}
\left\{
\begin{array}
[c]{l}%
x_{1,n+1}=f(ax_{2,n}+bx_{3,n})\\
x_{2,n+1}=cx_{1,n}+g(x_{3,n})\\
x_{3,n+1}=\alpha x_{1,n}+\beta
\end{array}
\right.  \label{sp3}%
\end{equation}
where $f,g:\mathbb{R\rightarrow R}$ are any non-constant functions and
$a,b,c,\alpha,\beta\in\mathbb{R}$ with $a,\alpha\not =0.$ We do not assume the
existence of partial inversions to solve either the first or the second
equation of the above system for $x_{2,n}$ or $x_{3,n}.$ We calculate%
\begin{align}
x_{1,n+3}  &  =f(ax_{2,n+2}+bx_{3,n+2})\nonumber\\
&  =f(acx_{1,n+1}+ag(x_{3,n+1})+b\alpha x_{1,n+1}+b\beta)\nonumber\\
&  =f((ac+b\alpha)x_{1,n+1}+ag(\alpha x_{1,n}+\beta)+b\beta) \label{sp3s}%
\end{align}

Equation (\ref{sp3s}) is derived from the system (\ref{sp3}) without any
partial inversions. In particular, if%
\begin{equation}
ac+b\alpha=0 \label{sco}%
\end{equation}
then (\ref{sp3s}) reduces to%
\begin{equation}
x_{1,n+3}=f(ag(\alpha x_{1,n}+\beta)+b\beta) \label{sp3s1}%
\end{equation}

Note that for $i=0,1,2$
\[
x_{1,3(k+1)+i}=f(ag(\alpha x_{1,3k+i}+\beta)+b\beta)
\]
so if $\{t_{k}\}$ is a solution of the first-order equation%
\begin{equation}
t_{k+1}=f(ag(\alpha t_{k}+\beta)+b\beta) \label{sp3o1}%
\end{equation}
with some initial value $t_{0}$ then
\[
t_{0}=x_{1,i}\Rightarrow t_{k}=x_{1,3k+i}\quad i=0,1,2
\]

Therefore, every solution of (\ref{sp3s1}) is obtained from 3 solutions of
(\ref{sp3o1}). In this sense, system (\ref{sp3}) may be reduced to a
first-order equation if (\ref{sco}) holds. If $\{x_{1,n}\}$ is a solution of
(\ref{sp3s1}) then for $n=1,2,3,\ldots$ the remaining components are found
passively as%
\[
x_{3,n}=\alpha x_{1,n-1}+\beta,\quad x_{2,n}=cx_{1,n-1}+g(x_{3,n-1}).
\]

\begin{remark}
\label{tok}The ideas discussed in this section do NOT apply to systems of
differential equations of type%
\begin{align}
x_{i}^{\prime}(t) &  =f_{i}(t,x_{1}(t),x_{i+1}(t),x_{i+2}(t),\ldots
,x_{k}(t)),\quad i=1,2,\ldots,k-1\label{odek}\\
x_{k}^{\prime}(t) &  =f_{k}(t,x_{1}(t))\nonumber
\end{align}
because the interdependence of variables in such systems is not the same as
that in (\ref{ske}). A more suitable analog of (\ref{odek}) is the system of
difference equations%
\begin{align}
\Delta x_{i,n} &  =f_{i}(n,x_{1,n},x_{i+1,n},x_{i+2,n},\ldots,x_{k,n}),\quad
i=1,2,\ldots,k-1\label{deltak}\\
\Delta x_{k,n} &  =f_{k}(n,x_{1,n})\nonumber
\end{align}
where $\Delta x_{i,n}=x_{i,n+1}-x_{i,n}.$ In this system the value of
$x_{i,n+1}$ depends on $x_{i,n}$ which is expressly not the case in
(\ref{ske}). Both of the systems (\ref{odek}) and (\ref{deltak}) require
partial inversions and are more difficult to fold into equations than
(\ref{ske}).

In general, difference and differential systems do not benefit equally from
folding because if a singularity occurs in the folding of a differential
system then it is less avoidable than the analogous situation in the discrete
case. Simply put, flows of differential systems are continuous and may thus
cross over or into singularity manifolds of the folding rather than jumping
over them. For instance, when Lorenz's system is folded (see \cite{Linz}) the
third order jerk function has a singularity in the terms $x^{\prime}/x$ that
the familiar, butterfly-shaped flows cannot seem to avoid. On the other hand,
the occurrence of singularities in the folding is not always fatal if the main
flows are not affected. For instance, the Volterra predator-prey model%
\[
\left\{
\begin{array}
[c]{l}%
x^{\prime}=x(a-by)\\
y^{\prime}=y(c-dx)
\end{array}
\right.
\]
where $a,b,c,d>0$ folds as follows: $by=a-x^{\prime}/x$ so%
\[
x^{\prime\prime}=ax^{\prime}-bx^{\prime}y-bxy^{\prime}=\frac{(x^{\prime})^{2}%
}{x}+(x^{\prime}+ax)(c-dx)
\]

Although a singularity exists in the folding on the y-axis this does not
involve the positive solutions that are important in modelling. Existence and
behavior of solutions for the above second-order equation may be tied to
similar issues about the flows of the Volterra system; such issues require
closer scrutiny in the future.
\end{remark}

%\bigskip

\end{document}